\documentclass[12pt]{article}
\input epsf.tex

\usepackage{amsmath}
\usepackage{amsthm}
\usepackage{amsfonts}
\usepackage{amssymb}
\usepackage{graphicx}
\usepackage{latexsym}

\usepackage{amsmath,amsthm,amsfonts,amssymb}

\usepackage{epsfig}

\theoremstyle{plain}
\newtheorem{thm}{Theorem}[section]
\newtheorem{prop}[thm]{Proposition}
\newtheorem{lem}[thm]{Lemma}
\newtheorem{cor}[thm]{Corollary}

\theoremstyle{definition}
\newtheorem{defn}{Definition}
\theoremstyle{remark}
\newtheorem{remark}{Remark}
\newtheorem{example}{Example}

\advance \topmargin by -\headheight
\advance \topmargin by -\headsep
\evensidemargin \oddsidemargin
\def\cA{{\cal A}}
\def\bA{{\bar A}}
\def\balpha{{\bar \alpha}}

\def\bB{{\bar B}}
\def\bbeta{{\bar \beta}}

\def\bE{{\mathbb E}}

\def\sF{{\mathcal F}}
\def\FA{{\mathcal{AP}}}

\def\tphi{{\bar\phi}}

\def\cG{{\mathcal G}}

\def\N{{\mathbb N}}

\def\sP{{\mathcal P}}

\def\R{{\mathbb R}}

\def\sym{{\textrm{Sym}}}

\def\sT{{\mathcal T}}

\def\v{{\vec v}}
\def\Var{{\textrm{Var}}}

\def\chix{{\raise.5ex\hbox{$\chi$}}}
\def\bxi{{\bar \xi}}
\def\bX{{\bar X}}

\def\Z{{\mathbb Z}}

\begin{document}
\title{Measure conjugacy invariants for actions of countable sofic groups}
\author{Lewis Bowen\footnote{email:lpbowen@math.hawaii.edu} \\ University of Hawaii}
\begin{abstract}
Sofic groups were defined implicitly by Gromov in [Gr99] and explicitly by Weiss in [We00]. All residually finite groups (and hence all linear groups) are sofic. The purpose of this paper is to introduce, for every countable sofic group $G$, a family of measure-conjugacy invariants for measure-preserving $G$-actions on probability spaces. These invariants generalize Kolmogorov-Sinai entropy for actions of amenable groups. They are computed exactly for Bernoulli shifts over $G$, leading to a complete classification of Bernoulli systems up to measure-conjugacy for many groups including all countable linear groups. Recent rigidity results of Y. Kida and S. Popa are utilized to classify Bernoulli shifts over mapping class groups and property (T) groups up to orbit equivalence and von Neumann equivalence respectively.
\end{abstract}
\maketitle

\noindent

\section{Introduction}
This paper is motivated by an old and central problem in measurable dynamics: given two dynamical systems, determine whether or not they are measurably-conjugate, i.e., isomorphic. Let us set some notation.

A {\bf dynamical system} (or system for short) is a triple $(G,X,\mu)$ where $(X,\mu)$ is a probability space and $G$ is a group acting by measure-preserving transformations on $(X,\mu)$. We will also call this a {\bf dynamical system over $G$}, a {\bf $G$-system} or an {\bf action of $G$}. In this paper, $G$ will always be a discrete countable group. Two systems $(G,X,\mu)$ and $(G,Y,\nu)$ are {\bf isomorphic} (i.e., {\bf measurably conjugate}) if and only if there exist conull sets $X' \subset X, Y' \subset Y$ and a bijective measurable map $\phi:X'\to Y'$ such that $\phi^{-1}:Y'\to X'$ is measurable, $\phi_*\mu=\nu$ and $\phi(gx)=g\phi(x) \forall g\in G, x \in X'$. 

A special class of dynamical systems called {\bf Bernoulli systems} or {\bf Bernoulli shifts} has played a significant role in the development of the theory as a whole because it was the problem of trying to classify them that motivated Kolmogorov to introduce the mean entropy of a dynamical system over $\Z$ [Ko58, Ko59].  That is, Kolmogorov defined for every system $(\Z,X,\mu)$ a number $h(\Z,X,\mu)$ called the {\bf mean entropy} of $(\Z,X,\mu)$ that quantifies, in some sense, how ``random'' the system is. His definition was modified by Sinai [Si59]; the latter has become standard. 


Bernoulli shifts also play an important role in this paper, so let us define them. Let $(K,\kappa)$ be a standard Borel probability space. For a discrete countable group $G$, let $K^G = \prod_{g \in G} K$ be the set of all functions $x: G \to K$ with the product Borel structure and let $\kappa^G$ be the product measure on $K^G$. The group $G$ acts on $K^G$ by $(gx)(f)=x(g^{-1}f)$ for $x \in K^G$ and $g,f \in G$. This action is measure-preserving. The system $(G,K^G,\kappa^G)$ is the {\bf Bernoulli shift over $G$ with base $(K,\kappa)$}. It is nontrivial if $\kappa$ is not supported on a single point.

Before Kolmogorov's seminal work [Ko58, Ko59], it was unknown whether all nontrivial Bernoulli shifts over $\Z$ were measurably conjugate to each other. He proved that $h(\Z,K^\Z,\kappa^\Z)=H(\kappa)$ where $H(\kappa)$, the {\bf entropy of $\kappa$} is defined as follows. If there exists a finite or countably infinite set $K' \subset K$ such that $\kappa(K')=1$ then
$$H(\kappa)= - \sum_{k\in K'} \mu(\{k\}) \log\big( \mu\big(\{k\} \big) \big)$$
where we follow the convention $0\log(0)=0$. Otherwise, $H(\kappa)=+\infty$. Thus two Bernoulli shifts over $\Z$ with different base measure entropies cannot be measurably conjugate. 

The converse was proven by D. Ornstein in the groundbreaking papers [Or70a, Or70b]. That is, he proved that if two Bernoulli shifts $(\Z, K^\Z, \kappa^\Z), (\Z,L^\Z,\lambda^\Z)$ are such that $H(\kappa)=H(\lambda)$ then they are isomorphic.

 In [Ki75], Kieffer proved a Shannon-McMillan theorem for actions of a countable amenable group $G$. In particular, he extended the definition of mean entropy from $\Z$-systems to $G$-systems. It is then not difficult to show from that Kolmogorov's theorem extends to Bernoulli shifts over $G$. 

In the landmark paper [OW87], Ornstein and Weiss extended most of the classical entropy theory from $\Z$-systems to $G$-systems where $G$ is any countable amenable group. (This paper also contains many results for nondiscrete amenable groups). In particular, they proved that if two Bernoulli shifts $(G, K^G, \kappa^G)$, $(G,L^G,\lambda^G)$ over a countably infinite amenable group $G$ are such that $H(\kappa)=H(\lambda)$ then they are isomorphic. Thus Bernoulli shifts over $G$ are completely classified by base measure entropy.

Let us say that a group $G$ is {\bf Ornstein} if $H(\kappa)=H(\lambda)$ implies $(G, K^G, \kappa^G)$ is isomorphic to $(G,L^G,\lambda^G)$ where $(K,\kappa)$ and $(L,\lambda)$ are any two standard Borel probability spaces. By the above, all countably infinite amenable groups are Ornstein. Stepin proved that any countable group that contains an Ornstein subgroup is itself Ornstein [St75]. This paper is not widely available; but a proof is also supplied in [Bo08b]. It is apparently unknown whether or not every countably infinite group is Ornstein. But an open case is that of Ol'shanskii's monsters [Ol91].

At the end of [OW87], Ornstein and Weiss presented a curious example suggesting that there might not be a reasonable entropy theory for nonamenable groups. It pertains to a well-known fundamental property of entropy: it is nonincreasing under factor maps. To explain, let $(G,X,\mu)$ and $(G,Y,\nu)$ be two systems. A map $\phi:X \to Y$ is a factor if $\phi_*\mu=\nu$ and $\phi(gx)=g\phi(x)$ for a.e. $x\in X$ and every $g\in G$. If $G$ is amenable then the mean entropy of a factor is less than or equal to the mean entropy of the source. This is essentially due to Sinai [Si59]. So if $K_n=\{1,\dots,n\}$ and $\kappa_n$ is the uniform probability measure on $K_n$ then $(G,K_2^G,\kappa_2^G)$, which has entropy $\log(2)$, cannot factor onto $(G,K_4^G,\kappa_4^G)$, which has entropy $\log(4)$.

The argument above fails if $G$ is nonamenable. Indeed, let $G=\langle a,b \rangle$ be a rank 2 free group. Identify $K_2$ with the group $\Z/2\Z$ and $K_4$ with $\Z/2\Z \times \Z/2\Z$. Then 
$$\phi(x)(g) := \big(x(g)+x(ga), x(g) + x(gb) \big) ~\forall x \in K_2^G, g\in G$$
is a factor map from $(G,K_2^G,\kappa_2^G)$ onto $(G,K_4^G,\kappa_4^G)$. This is Ornstein-Weiss' example.

There is an obvious factor map from $(G,K_4^G,\kappa_4^G)$ onto $(G,K_2^G,\kappa_2^G)$ (for any group $G$), so the authors speculated that if $G=\langle a,b \rangle$ then these two Bernoulli shifts might be measurably conjugate. We now know that this is false. The paper [Bo08a] introduced an invariant for dynamical systems over a free group that behaves similarly to Kolmogorov-Sinai entropy. In particular, it shows that Bernoulli shifts over a free group are completely classified by base measure entropy. 

The purpose of this paper is define, for every countable sofic group $G$, a family of isomorphism invariants that enables us to completely classify Bernoulli shifts over a countable sofic Ornstein group. It is expected that, as with Kolmogorov-Sinai entropy, these invariants will have an impact broader than this initial application to Bernoulli shifts. So what is a sofic group? Before stating the definition, let us note a few facts. These groups were defined implicitly by Gromov in [Gr99] and explicitly by Weiss in [We00]. An almost immediate consequence of the definition is that all residually amenable groups are sofic. In particular, since linear groups (i.e., subgroups of $GL_n(F)$ where $F$ is a field) are residually finite (by [Ma40]) they are sofic. It is unknown whether every countable group is sofic but an unresolved case is that of the universal Burnside group on a finite set of generators. Pestov has written a beautiful up-to-date survey [Pe08] on sofic groups and their siblings, hyperlinear groups.

\begin{defn}\label{defn:sofic}
Let $G$ be a countable group. For $m \ge 1$, let $\sym(m)$ denote the full symmetric group on $\{1,\dots, m\}$. Let $\sigma:G \to \sym(m)$ be a map. $\sigma$ is not assumed to be a homomorphism! For $F \subset G$, let $V(F) \subset \{1,\dots,m\}$ be the set of all elements $v$ such that for all $f_1,f_2 \in F$, 
$$\sigma(f_1)\sigma(f_2)v = \sigma(f_1f_2)v$$
and $\sigma(f_1)v \ne \sigma(f_2)v$ if $f_1 \ne f_2$. $\sigma$ is an {\bf $(F,\epsilon)$-approximation} to $G$ if $|V(F)| \ge (1-\epsilon)m$. 

Let $\Sigma=\{\sigma_i\}_{i=1}^\infty$ be a sequence of maps $\sigma_i:G \to \sym(m_i)$. Then $\Sigma$ is a {\bf sofic approximation} to $G$ if each $\sigma_i$ is an $(F_i,\epsilon_i)$-approximation to $G$ for some $(F_i,\epsilon_i)$ where $F_i \subset F_{i+1}$ for all $i$, $\cup_i F_i = G$ and $\epsilon_i \to 0$ as $i\to \infty$. $G$ is {\bf sofic} if there exists a sofic approximation to $G$. 
\end{defn}

\begin{example}
If $G$ is residually finite then there exists a sequence $\{N_i\}$ of finite-index normal subgroups of $G$ with $N_{i+1}<N_i$ for all $i$ and $\cap_i N_i = \{e\}$. Let $\sigma_i: G \to \sym(G/N_i)$ be the canonical homomorphism given by the action of $G$ on $G/N_i$. Then $\{\sigma_i\}$ is a sofic approximation to $G$.
\end{example}

\begin{example}
If $G$ is amenable then there exists an increasing sequence $\{F_i\}$ of finite subsets of $G$ such that $\bigcup_i F_i = G$ and for every finite $K \subset G$
$$\lim_{i\to\infty} \frac{|KF_i \Delta F_i|}{|F_i|} = 1.$$
Let $\sigma_i:G \to \sym(F_i)$ be any map such that if $f\in F_i$, $g\in G$ and $gf \in F_i$ then $\sigma_i(g)f=gf$. Then $\{\sigma_i\}$ is a sofic approximation to $G$.
\end{example}

In section \ref{sec:invariants}, we define the entropy of a system $(G,X,\mu)$ with respect to a sofic approximation $\Sigma$. It is denoted $h(\Sigma,G,X,\mu)$. The proof that this entropy is invariant under measure-conjugacy occupies sections \ref{sec:spaceofpartitions} - \ref{sec:monotone}. In section \ref{sec:Bernoulli}, it is proven that $h(\Sigma,G,K^G,\kappa^G)=H(\kappa)$ whenever $H(\kappa)<\infty$. This implies the next result.

\begin{thm}\label{thm:mainapp}
Let $G$ be a countable sofic group and let $(K_1,\kappa_1), (K_2,\kappa_2)$ be standard Borel probability spaces such that $H(\kappa_1)+H(\kappa_2) <\infty$. If $(G, K_1^G, \kappa_1^G)$ is isomorphic to $(G, K_2^G, \kappa_2^G)$ then $H(\kappa_1)=H(\kappa_2)$.
\end{thm}

In section \ref{sec:Bernoulli}, it is shown that if $G$ is also Ornstein then the finiteness condition in the above theorem can be removed. Thus:

\begin{cor}\label{cor:stepin}
Let $G$ be a countable sofic Ornstein group. Let $(K_1, \kappa_1), (K_2, \kappa_2)$ be standard Borel probability spaces. Then $(G, K_1^G, \kappa_1^G)$ is isomorphic to $(G, K_2^G,\kappa_2^G)$ if and only if $H(\kappa_1)=H(\kappa_2)$. 
\end{cor}

To make a contrast, recall that two systems $(G,X,\mu)$, $(G,Y,\nu)$ are {\bf weakly isomorphic} if $(G,X,\mu)$ is a factor of $(G,Y,\nu)$ and $(G,Y,\nu)$ is a factor of $(G,X,\mu)$. The next theorem is proven in [Bo08b].
\begin{thm}\label{thm:weaklyisomorphic}
Let $G$ be a countable group that contains a nonabelian free subgroup. Let $(K_1,\kappa_1), (K_2,\kappa_2)$ be any two nontrivial standard Borel probability spaces. Then $(G, K_1^G, \kappa_1^G)$ is weakly isomorphic to $(G, K_2^G, \kappa_2^G)$.
\end{thm}
The main ingredient in the proof is Ornstein-Weiss' example. In section \ref{sec:Bernoulli}, this is used to prove:
\begin{thm}\label{thm:generating}
Let $G$ be a countable sofic group that contains a nonabelian free subgroup. Let $(K,\kappa)$ be a standard Borel probability space with $H(\kappa)=+\infty$. Then $(G, K^G,\kappa^G)$ does not have a finite-entropy generating partition.
\end{thm}
The conclusion to this theorem is well-known to hold if $G$ is amenable. It is apparently unknown whether this result holds for all countable groups.
 
Let us consider the special case in which $G$ is a countably infinite linear group. Then every finitely generated subgroup of $G$ is residually finite by [Ma40]. Hence, $G$ is sofic. By the celebrated Tits alternative [Ti72], any finitely generated subgroup of $G$ is either virtually solvable (and hence amenable) or contains a nonabelian free group. Thus either $G$ contains a nonabelian free subgroup or it is amenable. In either case, it is Ornstein. Therefore, the following is proven.
\begin{cor}
Let $G$ be a countably infinite linear group. If $(K_1,\kappa_1), (K_2,\kappa_2)$ are standard Borel probability spaces then $(G, K_1^G,\kappa_1^G)$ is isomorphic to $(G, K_2^G,\kappa_2^G)$ if and only if $H(\kappa_1)=H(\kappa_2)$. $G$ is nonamenable if and only if every two nontrivial Bernoulli shifts over $G$ are weakly isomorphic. If $H(\kappa_1)=+\infty$ then there are no finite-entropy generating partitions for $(G, K_1^G,\kappa_1^G)$. 
\end{cor}

\subsection{Conjugation up to automorphisms}
There following definition is important in the applications that follow.

\begin{defn}
Two systems $(G_1, X_1,\mu_1)$ and $(G_2, X_2,\mu_2)$ are {\bf conjugate up to automorphisms} if there exists an isomorphism $\Phi:G_1 \to G_2$ and a measure-space isomorphism $\phi:X'_1 \to X'_2$ (where $X'_i$ is a conull subset of $X_i$) such that $\phi(gx)=\Phi(g)\phi(x)$ for every $g\in G_1$ and $x\in X'_1$.

For example, let $(G,X,\mu)$ be a system and let $a:G \to G$ be an automorphism. Let $(X^a,\mu^a)$ be a copy of $(X,\mu)$. Define an action of $G$ on $(X^a,\mu^a)$ by $g\cdot x = a(g)x$ for $g\in G, x\in X^a$. Then $(G,X,\mu)$ and $(G,X^a,\mu^a)$ are conjugate up to automorphisms. It is possible that they are not isomorphic as $G$-systems.
\end{defn}


\begin{thm}\label{thm:mainapp2}
Let $G$ be a countable sofic group and let $(K_1,\kappa_1), (K_2,\kappa_2)$ be standard Borel probability spaces such that $H(\kappa_1)+H(\kappa_2) <\infty$. If $(G, K_1^G, \kappa_1^G)$ is conjugate up to automorphisms to $(G, K_2^G, \kappa_2^G)$ then $H(\kappa_1)=H(\kappa_2)$. 
\end{thm}

\begin{cor}\label{cor:stepin2}
Let $G$ be a countable sofic Ornstein group. Let $(K_1, \kappa_1), (K_2, \kappa_2)$ be standard Borel probability spaces. Then $(G, K_1^G, \kappa_1^G)$ is conjugate up to automorphisms to $(G, K_2^G,\kappa_2^G)$ if and only if $H(\kappa_1)=H(\kappa_2)$. 
\end{cor}

This theorem and its corollary follow from theorem \ref{thm:K}, proposition \ref{prop:Bernoulli} and lemma \ref{lem:auto} below.

\subsection{Orbit equivalence and von Neumann equivalence}

The purpose of this subsection is to show that if the group $G$ satisfies certain additional hypotheses then the results above can be used to classify Bernoulli shifts over $G$ up to orbit equivalence and even up to von Neumann equivalence. To recall the definitions, for $i=1,2$ let $(G_i,X_i,\mu_i)$ be a system.


\begin{defn}
$(G_1, X_1,\mu_1)$ and $(G_2, X_2,\mu_2)$ are {\bf orbit-equivalent (OE)} if there exist conull sets $X'_1 \subset X_1$, $X'_2\subset X_2$ and a measure-space isomorphism $\phi:X'_1\to X'_2$ such that for all $x\in X'_1$, $\phi(G_1x)=G_2\phi(x)$.
\end{defn}

Orbit equivalence was introduced implicitly in [Si55] and explicitly in [Dy59] where it was shown that all free ergodic actions of $\Z$ are orbit equivalent. This was extended to countable amenable groups in [OW80, CFW81]. In the last decade, a number of striking rigidity results in orbit equivalence theory have been proven. These imply that under special additional hypotheses, orbit equivalence implies conjugacy up to automorphisms. For example, S. Popa proved [corollary 1.3, Po08] that if $G_1, G_2$ are two countably infinite groups, $G_1$ is nonamenable and $G_1 \times G_2$ has no nontrivial finite normal subgroups then any two Bernoulli shifts over $G_1 \times G_2$ are orbit equivalent if and only if they are conjugate up to automorphisms. Thus theorem \ref{thm:mainapp2} and corollary \ref{cor:stepin2} imply
\begin{cor}\label{cor:1}
Let $G_1, G_2$ be countably infinite groups and suppose $G_1$ is nonamenable. Suppose that $G:=G_1 \times G_2$ is sofic and has no nontrivial finite normal subgroups. Let $(K_1, \kappa_1), (K_2, \kappa_2)$ be standard Borel probability spaces with $H(\kappa_1)+H(\kappa_2)<\infty$. If $(G, K_1^G, \kappa_1^G)$ is orbit-equivalent to $(G, K_2^G,\kappa_2^G)$ then $H(\kappa_1)=H(\kappa_2)$. 

If $G$ is also Ornstein, then $H(\kappa_i)$ is allowed to be infinite and the converse also holds. That is: $(G,K_1^G,\kappa_1^G)$ is OE to $(G,K_2^G,\kappa_2^G)$ if and only if $H(\kappa_1)=H(\kappa_2)$.
\end{cor}

Y. Kida proved [Ki08, theorem 1.4] that if $G$ is the mapping class group of a genus $g$, $n$-holed surface $S_{g,n}$ for some $(g,n)$ with $3g+n-4>0$ and $(g,n) \notin \{(1,2), (2,0)\}$, then any two Bernoulli shifts over $G$ are orbit equivalent if and only if they are conjugate up to automorphisms. By [Gr74, see also Iv86], mapping class groups are residually finite and hence, sofic. It is well-known that they contain infinite cyclic subgroups and hence, are Ornstein. So corollary \ref{cor:stepin2} implies:
\begin{cor}
Let $G$ be the mapping class group of a genus $g$, $n$-holed surface $S_{g,n}$ for some $(g,n)$ with $3g+n-4>0$ and $(g,n) \notin \{(1,2), (2,0)\}$. Let $(K_1, \kappa_1), (K_2, \kappa_2)$ be standard Borel probability spaces. Then $(G, K_1^G, \kappa_1^G)$ is orbit-equivalent to $(G, K_2^G,\kappa_2^G)$ if and only if $H(\kappa_1)=H(\kappa_2)$. 
\end{cor}
It is worth pointing out that Popa's and Kida's results are much more general than we have used here.


Let us now turn our attention towards von Neumann equivalence. A system $(G,X,\mu)$ gives rise in a natural way to a von Neumann algebra $L^\infty(X,\mu) \rtimes G$ called the {\bf group measure space} or {\bf crossed product construction} of Murray and von Neumann [MvN36]. If $G$ is infinite and the action is free and ergodic then $L^\infty(X,\mu) \rtimes G$ is a $II_1$ factor, a highly noncommutative infinite-dimensional algebra with a positive trace. It is a fundamental problem in the theory of von Neumann algebras to classify type $II_1$-factors up to isomorphism in terms of the group/action data. This motivates the next definition.
\begin{defn}
Two systems $(G_1,X_1,\mu_1)$ and $(G_2,X_2,\mu_2)$ are {\bf von Neumann equivalent (vNE)} if $L^\infty(X_1,\mu_1) \rtimes G_1$ is isomorphic to $L^\infty(X_2,\mu_2) \rtimes G_2$.
\end{defn}
It was shown in [Si55] that orbit equivalence implies von Neumann equivalence, indeed this insightful discovery motivated the study of orbit equivalence. In [Co76], A. Connes proved that all $II_1$ factors arising from actions of amenable groups are isomorphic. By contrast, nonamenable groups were used to produce large families of nonisomorphic factors in [MvN43], [Dy63], [Sc63], [Mc69], [Co75]. 

In a series of groundbreaking papers [Po06-Po08], S. Popa established a variety of vNE rigidity results. These posit that under certain additional hypotheses von Neumann equivalence implies conjugacy up to automorphisms. The survey [Po07] covers many of these developments. For example, in [Po06, corollary 0.2], it is proven that if $G$ is a countably infinite property (T) group such that every nontrivial conjugacy class is infinite (this is abbreviated as {\bf ICC}) then two Bernoulli shifts over $G$ are von Neumann equivalent if and only if they are conjugate up to automorphisms. Thus theorem \ref{thm:mainapp2} and corollary \ref{cor:stepin2} imply the following.

\begin{cor}
Let $G$ be a countably infinite ICC sofic property (T) group. Let $(K_1, \kappa_1)$, $(K_2, \kappa_2)$ be standard Borel probability spaces with $H(\kappa_1)+H(\kappa_2)<\infty$. If $(G, K_1^G, \kappa_1^G)$ is von Neumann equivalent to $(G, K_2^G,\kappa_2^G)$ then $H(\kappa_1)=H(\kappa_2)$. 

If, in addition, $G$ is Ornstein then $H(\kappa_1)$ and $H(\kappa_2)$ are allowed to be infinite and the converse also holds. That is, $L^\infty(K_1^G,\kappa_1^G) \rtimes G$ is isomorphic to $L^\infty(K_2^G,\kappa_2^G) \rtimes G$ if and only if $H(\kappa_1)=H(\kappa_2)$.
\end{cor}

For example, if $G=PSL_n(\Z)$ for $n>2$ then $G$ is a countably infinite, ICC, sofic, property (T), Ornstein group. The above result is a special case of a more general theorem. To state it we will need some definitions that are not widely known. For fundamental results related to these definitions, the reader is referred to [Po06].


\begin{defn}
A group $G$ is {\bf $w$-rigid} if it contains an infinite normal subgroup with the relative property (T) of Kazhdan-Margulis (in other words, $(G,H)$ is a property (T)-pair, see [Ma82], [dHv89]). For example, all infinite groups with property (T) are $w$-rigid.
\end{defn}

\begin{defn}
A subgroup $H<G$ is {\bf $wq$-normal} if for every intermediate subgroup $H<H' < G$ with $H' \ne G$, there exists an element $g\in G$ such that $|gH'g \cap H'|=+\infty$. $wq$-normal stands for ``weakly-quasinormal''.
\end{defn}

\begin{defn}
A group $G$ is in the class $w\sT_0$ if it contains a subgroup $H$ such that
\begin{itemize}
\item $(G,H)$ is a property (T) pair,
\item $H$ is not virtually abelian and
\item $H$ is $wq$-normal in $G$.
\end{itemize}
\end{defn}

\begin{cor}
Let $G$ be an ICC sofic group. Suppose one of the following conditions hold:
\begin{enumerate}
\item $G$ is $w$-rigid or in the class $w\sT_0$.
\item There is a nonamenable subgroup $H<G$ such that $C(H)$, the centralizer of $H$, is $wq$-normal in $G$ and is not virtually abelian. 
\end{enumerate}
Let $(K_1,\kappa_1), (K_2,\kappa_2)$ be two standard Borel probability spaces such that $H(\kappa_1)+H(\kappa_2)<\infty$. If $(G,K_1^G,\kappa_1^G)$ is von Neumann equivalent to $(G,K_2^G,\kappa_2^G)$ then $H(\kappa_1)=H(\kappa_2)$. 

If, in addition, $G$ is Ornstein (for example, if $G$ is linear) then $H(\kappa_1)$ and $H(\kappa_2)$ are allowed to be infinite and the converse also holds. That is, $L^\infty(K_1^G,\kappa_1^G) \rtimes G$ is isomorphic to $L^\infty(K_2^G,\kappa_2^G) \rtimes G$ if and only if $H(\kappa_1)=H(\kappa_2)$.

\end{cor}

\begin{proof}
This follows from theorem \ref{thm:mainapp2}, corollary \ref{cor:stepin2}, [Po06, corollary 0.2] (if condition (1.) holds), and [Po08, theorem 1.5] (if condition (2.) holds). 
\end{proof}


The following generalizes corollary \ref{cor:1}.

\begin{cor}
Let $G$ be a sofic group satisfying: $G$ has no nontrivial finite normal subgroups, $G$ contains infinite commuting subgroups $H, H'$ with $H$ nonamenable and $H'<G$ is $wq$-normal. Let $(K_1,\kappa_1), (K_2,\kappa_2)$ be two standard Borel probability spaces such that $H(\kappa_1)+H(\kappa_2)<\infty$. If $(G, K_1^G,\kappa_1^G)$ is orbit-equivalent to $(G, K_2^G,\kappa_2^G)$ then $H(\kappa_1)=H(\kappa_2)$.

If $G$ is also Ornstein then $H(\kappa_1)$ and $H(\kappa_2)$ are allowed to be infinite and the converse also holds. 
\end{cor}

\begin{proof}
This follows from theorem \ref{thm:mainapp2}, corollary \ref{cor:stepin2} and [Po08, corollary 1.3].
\end{proof}

\section{The invariants}\label{sec:invariants}
In this section, we define the new invariants and state the main theorem. So fix a countable group $G$. In this paper, all partitions $\alpha=(A_1,\dots)$ of a probability space $(X,\mu)$ are measurable and at most countable.

If $\Sigma=\{\sigma_i\}$ is a sofic approximation of $G$, $(G,X,\mu)$ is a system and $\alpha$ is a finite partition of $X$ then the $\Sigma$-entropy rate of $\alpha$ is, roughly speaking, the exponential rate of growth of the number of partitions $\beta$ on $\{1,\ldots,m_i\}$ that approximate $\alpha$. So we begin by making precise a notion of approximation for such partitions.

\begin{defn}\label{defn:FA}
Let $(G,X,\mu)$ be a system and $\alpha=(A_1, A_2, \dots )$ an ordered partition of $X$. Let $\sigma: G \to \sym(m)$ be a map, $\zeta$ be the uniform probability measure on $\{1,\dots,m\}$ and $\beta =(B_1, B_2, \dots)$ be a partition of $\{1,\dots,m\}$.  

Let $F \subset G$ be finite. Given a function $\phi:F \to \N$, let $A_\phi = \bigcap_{f \in F} fA_{\phi(f)}$ and $B_\phi = \bigcap_{f\in F} \sigma(f)B_{\phi(f)}$. Define
$$d_F(\alpha,\beta)=\sum_{\phi:F \to \N} \Big| \mu( A_\phi) - \zeta(B_\phi )\Big|.$$
The above definitions make sense even if $\alpha=(A_1,\dots, A_u)$ or $\beta=(B_1,\dots,B_v)$ are finite: just set $A_i=B_j=\emptyset$ for $i > u$ and $j>v$.


If $\alpha=(A_1,\dots,A_u)$ is finite then for $\epsilon>0$, let $\FA(\sigma,\alpha: F,\epsilon)$ be the set of all ordered partitions $\beta=(B_1,\dots,B_u)$ of $\{1,\dots, m\}$ with the same number of atoms as $\alpha$ such that $d_F(\alpha,\beta)\le \epsilon$. $\FA$ stands for {\bf approximating partitions}.
\end{defn}

\begin{defn}
A {\bf map sequence} for $G$ is a sequence $\Sigma=\{\sigma_i\}_{i=1}^\infty$ of maps $\sigma_i:G \to \sym(m_i)$ such that $m_i \to \infty$ as $i\to\infty$. 
\end{defn}

\begin{defn}
Let $(G,X,\mu)$ be a system and $\Sigma:=\{\sigma_i\}_{i=1}^\infty$ a map sequence for $G$. For every finite partition $\alpha$ of $X$, $\epsilon >0$ and finite set $F \subset G$, let
\begin{eqnarray*}
H(\Sigma, \alpha: F, \epsilon)&=&\limsup_{i\to\infty} \frac{1}{m_i} \log  |\FA(\sigma_i,\alpha: F, \epsilon)|\\
H(\Sigma,\alpha:F) &=& \lim_{\epsilon \to 0} H(\Sigma, \alpha: F, \epsilon).\\
h(\Sigma,\alpha) &=& \inf_{F \subset G} H(\Sigma, \alpha: F).
\end{eqnarray*}
A few words about the above definitions are in order. If $\epsilon_1 \ge \epsilon_2$ then $\FA(\sigma,\alpha:F,\epsilon_1) \supset \FA(\sigma,\alpha:F,\epsilon_2)$, so the limit defining $H(\Sigma,\alpha:F)$ exists and equals the infimum over all $\epsilon>0$.

The infimum defining $h(\Sigma,\alpha)$ is over all finite sets $F \subset G$. We call $h(\Sigma,\alpha)$ the {\bf mean $\Sigma$-entropy} of $\alpha$. Note that if $F_1 \subset F_2$ then $\FA(\sigma,\alpha:F_1,\epsilon) \supset \FA(\sigma,\alpha:F_2,\epsilon)$. Hence if $\{F_n\}$ is any sequence of finite subsets of $G$ with $F_n \subset F_{n+1}$ for all $n$ and $\cup_n F_n =G$ then $h(\Sigma,\alpha) = \lim_{n\to\infty} H(\Sigma,\alpha:F_n)$. 

It is possible that $\FA(\sigma_i,\alpha:F,\epsilon)$ is empty. In this case, we interpret $\log(0)=-\infty$. Thus, it is apriori possible that $h(\Sigma,\alpha)=-\infty$. 
\end{defn}

In order to handle the case when $\alpha$ is an infinite partition, we need to review some standard definitions.
\begin{defn}
Let $(X,\mu)$ be a probability space and let $\alpha=(A_1, A_2, \dots)$ be a measurable partition of $X$ into at most countably many sets (each of which is called an {\bf atom} of $\alpha$). The {\bf entropy} of $\alpha$ is
$$H(\alpha) = -\sum_{i=1}^\infty \mu(A_i)\log(\mu(A_i)).$$
\end{defn}

\begin{defn}
If $\alpha$ and $\beta$ are two partitions of $X$ then their {\bf join}, denoted $\alpha \vee \beta$, is defined by $\alpha \vee \beta = \{A \cap B~|~ A \in \alpha , B \in \beta\}$. 
\end{defn}

\begin{defn}\label{defn:chain}
If $\alpha, \beta$ are partitions of $X$ and for every $A \in \alpha$ there exists a $B \in \beta$ such that $\mu(A-B)=0$ (i.e. $A$ is a subset of $B$ up to a set of measure zero) then we say that $\alpha$ {\bf refines} $\beta$. Equivalently, $\beta$ is a {\bf coarsening} of $\alpha$. We denote this by $\alpha \ge \beta$. A {\bf chain} of $\alpha$ is a sequence $\{\alpha_n\}_{n=1}^\infty$ of finite partitions such that $\alpha_1 \le \alpha _2 \le \dots \le \alpha$ and $\bigvee_{i=1}^\infty \alpha_i = \alpha$. 

Often, we will abuse notation by writing $A \subset B$ to mean $\mu(A-B)=0$.
\end{defn}

\begin{defn}
Let $(G,X,\mu)$ be a system and $\Sigma:=\{\sigma_i\}_{i=1}^\infty$ a map sequence for $G$. Let $\alpha=(A_1,A_2,\dots)$ be a partition of $X$. For every $\epsilon >0$ and finite set $F \subset G$, let
\begin{eqnarray*}
H(\Sigma,\alpha: F) &=& \inf\Big\{ \lim_{n \to\infty} H(\Sigma, \alpha_n: F) ~:~ \{\alpha_n\}_{n=1}^\infty \textrm{ is a chain of } \alpha\Big\}.
\end{eqnarray*}
We will prove in section \ref{sec:uppersemi} that $H(\Sigma,\alpha:F) = \lim_{n\to\infty} H(\Sigma,\alpha_n:F)$ for any chain $\{\alpha_n\}$ of $\alpha$. As in the finite case, define
$h(\Sigma,\alpha) = \inf_{F \subset G} H(\Sigma, \alpha: F).$
\end{defn}

\begin{remark}
The above definitions admit two natural generalizations. First, the limsup in the definition of $H(\Sigma,\alpha:F,\epsilon)$ can be replaced with a liminf or an ultralimit. The latter is very natural from the perspective on sofic groups taken in [ES05] and [Pe08]. Second, each map $\sigma_i:G \to \sym(m_i)$ could be random. In this case, define
$$H(\Sigma,\alpha:F,\epsilon) = \limsup_{i\to\infty} \frac{1}{m_i} \log \bE\Big[ \big|\FA(\sigma_i,\alpha:F,\epsilon)\big|\Big]$$
where $\bE[\cdot]$ denotes expected value. This is used in [Bo09] to show that the $f$-invariant defined in [Bo08a] is a special case of $\Sigma$-entropy. All of the results in this paper remain true if these two generalizations are utilized (with only minor, obvious changes in the proofs). But for simplicity's sake, we will not make use of either generalization.
\end{remark}

In order to obtain a measure-conjugacy invariant, we need to focus on a special class of partitions described next.
\begin{defn}\label{defn:generating}
Let $(G,X,\mu)$ be a system and $\alpha$ a partition of $X$. Let $\Sigma_\alpha$ be the smallest $G$-invariant $\sigma$-algebra containing the atoms of $\alpha$. Then $\alpha$ is {\bf generating} if for any measurable set $A \subset X$ there exists a set $A' \in \Sigma_\alpha$ such that $\mu(A \Delta A')=0$. 
\end{defn}

The main theorem of this paper is:
\begin{thm}\label{thm:K}
Let $G$ be a countable sofic group.  Let $\Sigma=\{\sigma_i\}_{i=1}^\infty$ be a sofic approximation to $G$. Then if $(G,X,\mu)$ is any $G$-system and $\alpha, \beta$ are any two generating partitions of $X$ with $H(\alpha)+H(\beta)<\infty$ then $h(\Sigma, \alpha)=h(\Sigma,\beta)$. 
\end{thm}
This motivates the following definition:
\begin{defn}
Let $(G,X,\mu), \Sigma$ be as above. If $(G,X,\mu)$ has a generating partition $\alpha$ with $H(\alpha)<\infty$ then let $h(\Sigma, G,X,\mu)=h(\Sigma,\alpha)$. By the above theorem, $h(\Sigma,G,X,\mu)$ depends on the system $(G,X,\mu)$ only up to measure-conjugacy. 
\end{defn}


\begin{remark}
If there does not exist a generating partition $\alpha$ with $H(\alpha)<\infty$ then $h(\Sigma,G,X,\mu)$ is undefined. This differs from the classical case in which the mean entropy of a system is defined as the supremum of the mean entropy rate of $\alpha$ over all finite partitions $\alpha$. In general, it is possible that the supremum of $h(\Sigma,\alpha)$ over all finite partitions $\alpha$ is infinite even if $h(\Sigma,G,X,\mu)$ is finite. For example, theorem \ref{thm:weaklyisomorphic} implies that this occurs for the Bernoulli shift $(G,K^G,\kappa^G)$ if $G$ contains a nonabelian free group.

\end{remark}

\begin{remark}
It can be shown that if $G$ is amenable and if $\Sigma$ is any sofic approximation to $G$ then $h(\Sigma,G,X,\mu)$ is the classical mean entropy of $(G,X,\mu)$. Since we will not use this, we do not prove it.
\end{remark}

In section \ref{sec:Bernoulli} we show the following.
\begin{prop}\label{prop:Bernoulli}
Let $G$ be a countable sofic group. Let $\Sigma$ be a sofic approximation to $G$. Let $(K,\kappa)$ be a probability space with $H(\kappa)<\infty$. Then $h(\Sigma,G,K^G,\kappa^G) = H(\kappa)$.
\end{prop}
This proposition and theorem \ref{thm:K} imply theorem \ref{thm:mainapp}. Theorem \ref{thm:mainapp2} follows from this proposition, theorem \ref{thm:K} and the next lemma.

\begin{lem}\label{lem:auto}
Let $a:G \to G$ be an automorphism of a countable sofic group $G$. Let $\Sigma=\{\sigma_i\}$ be a sofic approximation and $(G,X,\mu)$ a $G$-system. Let $(X^a,\mu^a)$ be a copy of $(X,\mu)$. Define an action of $G$ on $(X^a,\mu^a)$ by $g\cdot x =a(g)x$. Then $h(\Sigma,G,X,\mu)=h(\Sigma^a,G,X^a,\mu^a)$ where $\Sigma^a$ is the sequence $\{\sigma_i \circ a\}_{i=1}^\infty$. 
\end{lem}
\begin{proof}

Let $\epsilon>0$, $F \subset G$ be finite. Let $\alpha$ be a finite partition of $X$. For $i\ge 0$, define $\FA(\sigma_i,\alpha:F,\epsilon)$ as in definition \ref{defn:FA}. 

Let $\alpha^a$ be the partition of $X^a$ corresponding to $\alpha$ (i.e., $\alpha^a$ is a copy of $\alpha$ in $X^a$). Let $\FA^a\big(\sigma_i\circ a,\alpha^a:a^{-1}(F),\epsilon\big)$ be defined as in definition \ref{defn:FA} - but with the system $(G,X^a,\mu^a)$ in place of $(G,X,\mu)$. An exercise in definition-chasing shows that $\FA(\sigma_i,\alpha:F,\epsilon)=\FA^a\big(\sigma_i\circ a,\alpha^a:a^{-1}(F),\epsilon\big)$. This implies $h(\Sigma,\alpha) = h(\Sigma^a,\alpha^a)$. Since this is true for every finite partition $\alpha$, the lemma follows.




\end{proof}

\subsection{An Outline of the Paper}

To prove theorem \ref{thm:K}, we first review standard definitions and results from classical entropy theory. This is section \ref{sec:standard}. Then, we develop a general theory for the space of partitions in section \ref{sec:spaceofpartitions}. This culminates in the establishment of simple criteria under which a function $f$ from the space of partitions to $\R \cup \{-\infty\}$ is guaranteed to be constant on the set of generating partitions. It is shown in sections \ref{sec:factor}-\ref{sec:monotone} that $h(\Sigma,\cdot)$ satisfies these conditions. This finishes the proof of theorem \ref{thm:K}. In section \ref{sec:Bernoulli}, we compute the entropy of a Bernoulli shift over a sofic group with respect to a sofic approximation and use it to prove proposition \ref{prop:Bernoulli}, theorem \ref{thm:mainapp}, corollary \ref{cor:stepin} and theorem \ref{thm:generating}. 

\section{Classical Entropy Theory in Brief}\label{sec:standard}
To prove the above results, we will need some basic facts from classical entropy theory. An expert could skip this section, referring back to it for notation if necessary. Fix a probability space $(X,\mu)$.


\begin{defn}\label{defn:conditional}
Let $\sF$ be a $\sigma$-algebra contained in the $\sigma$-algebra of all measurable subsets of $X$. Given a partition $\alpha$ of $X$, define the {\bf conditional information function} $I(\alpha|\sF):X \to \R$ by
$$I(\alpha|\sF)(x) = -\log\big(\mu(A_x|\sF)(x)\big)$$
where $A_x$ is the atom of $\alpha$ containing $x$. Here, if $A \subset X$ is measurable then $\mu(A|\sF):X \to \R$ is the conditional expectation of $\chi_{A}$, the characteristic function of $A$, with respect to the $\sigma$-algebra $\sF$. In other words, it is an $\sF$-measurable function such that for all $\sF$-measurable functions $f:X \to \R$,
$$\int_X \mu(A|\sF)(x)f(x) \, d\mu(x) = \int_X \chi_{A}(x)f(x) \, d\mu(x).$$
The {\bf conditional entropy of $\alpha$ with respect to $\sF$} is defined by
$$H(\alpha|\sF) = \int_X I(\alpha | \sF)(x) \, d\mu(x).$$
For simplicity, let $H(\alpha)=H(\alpha | \{X,\emptyset\})$ and $I(\alpha)=I(\alpha | \{X,\emptyset\})$. 

If $\beta$ is a partition then, by abuse of notation, we can identify $\beta$ with the $\sigma$-algebra equal to the set of all unions of partition elements of $\beta$. Through this identification, $I(\alpha|\beta)$ and $H(\alpha|\beta)$ are well-defined.
\end{defn}

\begin{lem}\label{lem:relative}
For any two partitions $\alpha, \beta$ and for any two $\sigma$-algebras $\sF_1, \sF_2$ with $\sF_1 \subset \sF_2$,
$$H(\alpha \vee \beta) = H(\alpha) + H(\beta|\alpha),$$
$$H(\alpha  | \sF_2) \le H(\alpha  | \sF_1)$$
with equality if and only if $\mu(A |\sF_2) = \mu(A |\sF_1)$ a.e. for every $A \in \alpha$. In particular $H(\alpha|\beta) \le H(\alpha)$ and equality occurs iff $\alpha$ and $\beta$ are {\bf independent} (i.e., $\forall A \in \alpha, B \in \beta, \mu(A\cap B)=\mu(A)\mu(B)$).
\end{lem}
\begin{proof}
This is well-known. For example, see [Gl03, Proposition 14.16, page 255].
\end{proof}

\section{The space of partitions}\label{sec:spaceofpartitions}

In order to prove theorem \ref{thm:K}, it is necessary to develop a general theory of the space of all partitions of a given probability space $(X,\mu)$ on which $G$ acts by measure-preserving transformations. In the case of finite partitions and finitely generated groups, the required results were proven in [Bo08a]. So fix a countable group $G$. Let $(G,X,\mu)$ be a $G$-system.

\begin{defn}
Two measurable partitions $\alpha, \beta$ of $X$ are {\bf equivalent} if for every $A \in \alpha$ there exists a $B\in \beta$ such that $\mu(A \Delta B)=0$. Let $\sP$ denote the set of equivalence classes of measurable partitions $\alpha$ of $X$ with $H(\alpha)<\infty$. Generally speaking, we will abuse notation and consider elements of $ \sP$ as partitions themselves. This is similar to the common abuse of considering an element $f \in L^p(X,\mu)$ to be a function when it is really an equivalence class of functions.
\end{defn}

\begin{defn}[Rohlin distance]\label{defn:rohlin}
Define $d:\sP \times \sP \to \R$ by
$$d(\alpha,\beta) = H(\alpha|\beta) + H(\beta|\alpha) = 2H(\alpha \vee \beta) - H(\alpha) - H(\beta).$$
By [Pa69, theorem 5.22, page 62] this defines a distance function.

$G$ acts isometrically on $\sP$ by $g\alpha=(gA_1,gA_2,\dots)$ where $g\in G$ and $\alpha=(A_1,A_2,\dots) \in \sP$.  I.e., if $g \in G$, $\alpha, \beta \in \sP$ then $d(g\alpha, g\beta) = d(\alpha,\beta)$.
\end{defn}


\begin{defn}\label{defn:cayley}
Let $S \subset G$. The {\bf left-Cayley graph} $\Gamma$ of $(G,S)$ is defined as follows. The vertex set of $\Gamma$ is $G$. For every $s \in S$ and every $g \in G$ there is a directed edge from $g$ to $sg$ labeled $s$. There are no other edges. $S$ {\bf generates $G$} if and only if $\Gamma$ is connected.

The {\bf induced subgraph} of a subset $F \subset G$ is the largest subgraph of $\Gamma$ with vertex set $F$. A subset $F \subset G$ is {\bf $S$-connected} if its induced subgraph in $\Gamma$ is connected (as an undirected graph).
\end{defn}

\begin{defn}
Let $\alpha$ and $\beta$ be partitions. If, for every atom $A \in \alpha$ there exists an atom $B \in \beta$ such that $\mu(A-B)=0$ (i.e., $A$ is contained in $B$ up to a set of measure zero) then we say $\alpha$ {\bf refines} $\beta$. Equivalently, $\beta$ is a {\bf coarsening} of $\alpha$. This is denoted by $\beta \le \alpha$. We will often abuse notation by writing $A \subset B$ to mean $\mu(A-B)=0$.
\end{defn}

\begin{defn}
If $F \subset G$ is finite, and $\alpha \in \sP$, let
$$\alpha^F = \bigvee_{g \in F} g\alpha.$$
Partitions $\alpha, \beta \in \sP$ are {\bf $S$-equivalent} if there exists finite $S$-connected sets $F_1, F_2 \subset G$ such that $e \in F_1 \cap F_2$, $\alpha \le \beta^{F_1}$ and $\beta \le \alpha^{F_2}$. If $S, \alpha,\beta$ are all finite and $S$ generates $G$ then this is equivalent to the definition of combinatorial equivalence given in [Bo08a].

To check that this defines an equivalence relation, let $\alpha,\beta,\gamma \in \sP$ and let $A,B,C,D \subset G$ be finite $S$-connected sets containing the identity element such that $\alpha \le \beta^A$, $\beta \le \alpha^B$, $\beta \le \gamma^C$ and $\gamma\le \beta^D$. Then
$$\alpha \le \beta^A \le (\gamma^C)^A = \gamma^{AC}$$
where $AC=\{ac~|~a\in A, c\in C\}$. Similarly, $\gamma \le \alpha^{DB}$. An easy exercise shows that $AC$ and $DB$ both contain the identity and are $S$-connected. Thus $\alpha$ is $S$-equivalent to $\gamma$. This shows that $S$-equivalence is an equivalence relation on $\sP$.
\end{defn}

The main result about $S$-equivalence is:
\begin{thm}\label{thm:dense}
Let $S\subset G$ generate $G$. If $\alpha, \beta \in \sP$ are generating partitions of the system $(G,X,\mu)$, then for every $\epsilon>0$ there exists $\alpha' \in \sP$ such that $\alpha'$ is $S$-equivalent to $\alpha$ and $d(\alpha',\beta)<\epsilon$. I.e., the $S$-equivalence class of $\alpha$ is dense in the space of all generating partitions.
\end{thm}
We will first need some lemmas.

\begin{lem}\label{lem:alpha^n}
If $F \subset G$ is finite then the function $\alpha \mapsto \alpha^F$ is continuous on $\sP$.
\end{lem}

\begin{proof}
Let $\sP^F = \prod_{f \in F} \sP$ be the product space. Define $\phi: \sP \to \sP^F$ by $\phi(\alpha) = (f\alpha)_{f \in F}$. Define $\psi:\sP^F \to \sP$ by $\psi\big( (\alpha_f)_{f \in F}\big) = \bigvee_{f \in F} \alpha_f$. Then the function $\alpha \mapsto \alpha^F$ is the composition of $\phi$ and $\psi$.

For fixed $g\in G$, the map $\alpha \mapsto g\alpha$ is continuous on $\sP$. So $\phi$ is continuous (where $\sP^F$ has the product topology). It is easy to see that $\psi$ is also continuous. For example, if $\alpha, \alpha', \beta,\beta' \in \sP$ then
\begin{eqnarray*}
d(\alpha \vee \beta,\alpha' \vee \beta')&=&H(\alpha \vee \beta|\alpha'\vee\beta') + H(\alpha'\vee\beta'|\alpha\vee\beta)\\
&\le& H(\alpha |\alpha'\vee\beta')+H(\beta|\alpha'\vee\beta') + H(\alpha'|\alpha\vee\beta)+H(\beta'|\alpha\vee\beta)\\
&\le&H(\alpha|\alpha')+H(\beta|\beta')+H(\alpha'|\alpha)+H(\beta'|\beta)\\
&=&d(\alpha,\alpha')+d(\beta,\beta').
\end{eqnarray*}
Similarly, if $(\alpha_f)_{f \in F}, (\beta_f)_{f\in F} \in \sP^F$ then $d(\bigvee_{f \in F} \alpha_f, \bigvee_{f \in F} \beta_f) \le \sum_{f \in F} d(\alpha_f,\beta_f)$. This proves that $\alpha \mapsto \alpha^F$ is continuous.
\end{proof}


\begin{lem}\label{lem:predense}
Let $\alpha \in \sP$ be a generating partition. Let $\beta \in \sP$ and $\epsilon > 0$. Let $S\subset G$ generate $G$. Then there exists a finite $S$-connected set $F \subset G$ and a partition $\gamma \in \sP$ such that $\gamma \le \alpha^F$ and $d(\beta,\gamma) \le \epsilon$.
\end{lem}

\begin{proof}
This is an easy exercise left to the reader. 
\end{proof}


\begin{lem}
Let $\alpha \in \sP$. Then for every $\epsilon>0$ there exists a $\delta>0$ such that if $\omega=\{X_L,X_S\}$ is a 2-atom partition of $X$ with $\mu(X_S)<\delta$ and $\xi$ is the partition
$$\xi=\{X_L\} \cup \{X_S \cap A~|~A \in \alpha\},$$
then $H(\xi) <\epsilon$.
\end{lem}

\begin{proof}
Let $\epsilon>0$. Let $\delta>0$ be small enough so that if $Y\subset X$ is any set with $\mu(Y)<\delta$ and if $\omega=\{X-Y,Y\}$ then
$$2H(\omega) + \int_Y I(\alpha) ~ d\mu < \epsilon$$
where $I(\alpha)$ is the information function of $\alpha$.

Now let $\omega=\{X_L,X_S\}$ with $\delta> \mu(X_S)$ and let $\xi$ be as above. Consider the information function $I(\xi)$. Note that $I(\xi)(x)=I(\omega)(x)$ if $x \in X_L$ and $I(\xi)(x)=I(\alpha \vee \omega)(x)$ if $x \in X_S$. It follows from the definition of the information function that $I(\alpha \vee \omega)= I(\alpha) + I(\omega|\alpha)$. Hence
\begin{eqnarray*}
H(\xi) &=& \int_{X_L} I(\omega) ~ d\mu + \int_{X_S} I(\alpha \vee \omega) ~ d\mu\\
&\le & H(\omega) + \int_{X_S} I(\alpha) + I(\omega|\alpha) ~d\mu\\
&\le& 2H(\omega) + \int_{X_S} I(\alpha) ~d\mu < \epsilon.
\end{eqnarray*}
\end{proof}

We can now prove theorem \ref{thm:dense}.

\begin{proof}[Proof of theorem \ref{thm:dense}]
Let $1>\epsilon>0$. Let $\delta>0$ be as in the previous lemma. By choosing $\delta$ smaller if necessary, we may assume $e^{-\delta}> 1/2$ and $\epsilon>\delta$.

Because $S$ generates $G$, the Cayley graph $\Gamma$ of $(G,S)$ is connected. Since $\beta$ is generating, there exists a finite $S$-connected set $F \subset G$ such that $e\in F$ and $H(\alpha|\beta^F)<\delta^2/2$. Let $\delta_2$ be a number with $\delta^2/2>\delta_2>0$. By lemma \ref{lem:predense} there exists a finite $S$-connected set $K \subset G$ and a partition $\gamma \in \sP$ such that $e\in K$, $\gamma \le \alpha^K$ and $d(\gamma,\beta) \le \delta_2$. By lemma \ref{lem:alpha^n}, the function $\alpha \mapsto \alpha^F$ is continuous. Hence, by choosing $\delta_2$ smaller if necessary, we may assume that $d(\gamma^F, \beta^F) \le \delta^2/2$. So,
\begin{eqnarray*}
H(\alpha|\gamma^F) &=& H(\alpha \vee \gamma^F) - H(\gamma^F)\\
&\le & H(\alpha |\gamma^F \vee \beta^F) + H(\gamma^F \vee \beta^F) - H(\gamma^F)\\
&\le & H(\alpha|\beta^F) +H(\beta^F|\gamma^F)\\
&\le & H(\alpha|\beta^F) + d(\beta^F,\gamma^F) \le \delta^2.
\end{eqnarray*}
Let
$$X_L =\{x\in X ~|~ I(\alpha|\gamma^F)(x) \le \delta\}.$$
Let $X_S=X-X_L$. Since 
\begin{eqnarray*}
\delta^2 &\ge& H(\alpha|\gamma^F) = \int I(\alpha|\gamma^F)(x) ~ d\mu(x) \ge \mu(X_S)\delta,
\end{eqnarray*}
$\mu(X_S) \le \delta$. Let $\omega=\{X_L, X_S\}$, $\xi=\{X_L\} \cup \{X_S \cap A~|~ A \in \alpha\}$ and $\lambda = \gamma \vee \xi$. We will show that $\lambda$ is $S$-equivalent to $\alpha$ and $d(\lambda, \beta)<2\epsilon$. 

Since the information function $I(\alpha|\gamma^F)$ is constant on each atom of $\alpha \vee \gamma^F$, $X_L$ is a union of atoms of $\alpha \vee \gamma^F \le (\alpha^K)^F=\alpha^{FK}$ where $FK$ is the product $FK=\{fk~:~f\in F, k\in K\}$. Note that $FK$ is $S$-connected and $e\in FK$. So $\omega=\{X_L,X_S\} \le \alpha^{FK}$. This implies $\xi \le \alpha \vee \omega \le \alpha^{FK}$. Hence $\lambda =\gamma\vee \xi\le \alpha^{FK}$.

On the other hand, $\alpha \le \gamma^F \vee \xi$. To see this, let $A \in \alpha$. Then
$$A = (A \cap X_L) \cup (A \cap X_S).$$
Since $A \cap X_S \in \xi$, it suffices to show that $A \cap X_L$ is a union of atoms in $\gamma^F \vee \xi$. By definition, $A \cap X_L$ is the union of $A \cap C $ over all atoms $C \in \gamma^F$ such that $\mu(A|C) \ge e^{-\delta} > 1/2$ where $\mu(A|C)=\frac{\mu(A\cap C)}{\mu(C)}$. But if $\mu(A|C) >1/2$ and if $A_2 \in \alpha$ also satisfies $\mu(A_2|C) \ge e^{-\delta}>1/2$ then $A_2=A$. Thus $A\cap C = X_L \cap C$. So $A \cap X_L $ is the union of $X_L\cap C$ over all atoms $C \in \gamma^F$ such that $\mu(A|C) \ge e^{-\delta}$. This proves that $A \cap X_L$ is the union of atoms in $\gamma^F \vee \xi$. Thus, $\alpha \le \gamma^F \vee \xi$ as claimed. Since $\lambda=\gamma\vee\xi$, this implies $\alpha \le \lambda^F$. Because $\lambda \le \alpha^{FK}$, $\alpha$ is $S$-equivalent to $\lambda$.

It remains to estimate $d(\lambda,\beta)$.
\begin{eqnarray*}
d(\lambda,\beta)&=&H(\lambda|\beta) + H(\beta|\lambda)\\
&\le& H(\gamma \vee\xi |\beta) + H(\beta|\gamma)\\
&\le&  H(\xi|\beta) + d(\gamma,\beta) \\
&\le& H(\xi) + \delta^2/2  \le \epsilon +\delta^2/2.
\end{eqnarray*}
The last inequality comes from the definition of $\delta$ at the beginning of this proof. Since $\delta<\epsilon$ and $\epsilon$ is arbitrary, this proves the theorem with $\alpha'=\lambda$.
\end{proof}

\subsection{Splittings}
In order to apply theorem \ref{thm:dense}, we will show that if $\alpha, \beta \in \sP$ are $S$-equivalent then they have a common ``$S$-splitting'' as defined next. 

\begin{defn}\label{defn:$S$-splitting}
Let $\alpha \in \sP$ be a partition. A {\bf simple $S$-splitting} of $\alpha$ is a partition $\sigma$ of the form $\sigma=\alpha \vee s\beta$ where $s\in S$ and $\beta$ is a coarsening of $\alpha$.

An {\bf $S$-splitting} of $\alpha$ is any partition $\sigma$ that can be obtained from $\alpha$ by a sequence of simple $S$-splittings. In other words, there exist partitions $\alpha_0, \alpha_1,\ldots,\alpha_m$ such that $\alpha_0=\alpha$, $\alpha_m = \sigma$ and $\alpha_{i+1}$ is a simple $S$-splitting of $\alpha_i$ for all $0 \le i < m$.
\end{defn}

\begin{lem}\label{lem:equivalent}
If $\beta$ is an $S$-splitting of $\alpha \in \sP$ then $\alpha$ is $S$-equivalent to $\beta$.
\end{lem}
\begin{proof}
It suffices to consider the special case in which $\beta$ is a simple $S$-splitting. Then $\alpha \le \beta \le \alpha \vee t\alpha=\alpha^{\{e,t\}}$ for some $t \in S$. Since $\{e,t\}$ is $S$-connected, this proves it.
\end{proof}

\begin{lem}\label{lem:$S$-splittings}
Let $S\subset G$. If $\alpha, \beta \in \sP$, $\alpha$ refines $\beta$ and $F \subset G$ is finite, $S$-connected and contains the identity element $e$ then $\alpha \vee \beta^F$ is an $S$-splitting of $\alpha$.
\end{lem}

\begin{proof}
We prove this by induction on $|F|$. If $|F|=1$ then $F=\{e\}$ and the statement is trivial. Let $f_0 \in F-\{e\}$ be such that $F_1=F-\{f_0\}$ is $S$-connected. To see that such an $f_0$ exists, choose a spanning tree for the induced subgraph of $F$. Let $f_0$ be any leaf of this tree that is not equal to $e$.

By induction, $\alpha_1 := \alpha \vee \beta^{F_1}$ is an $S$-splitting of $\alpha$. Since $F$ is $S$-connected, there exists an element $f_1 \in F_1$ and an element $s_1 \in S$ such that $s_1f_1=f_0$. Since $f_1 \in F_1$, $f_1\beta \le \alpha_1$. Thus
 $$\alpha \vee \beta^F = \alpha_1 \vee f_0\beta= \alpha_1 \vee s_1(f_1\beta)$$ is an $S$-splitting of $\alpha$.
\end{proof}

\begin{prop}\label{prop:$S$-splitting}
Let $S\subset G$. Let $\alpha,\beta \in \sP$ be $S$-equivalent. Let $\balpha, \bbeta$ be $S$-splittings of $\alpha, \beta$ respectively. Then there exists a partition $\gamma \in \sP$ that is an $S$-splitting of $\balpha$ and an $S$-splitting of $\bbeta$.
\end{prop}

\begin{proof}
By lemma \ref{lem:equivalent}, $\balpha$ and $\bbeta$ are $S$-equivalent. So it suffices to prove that there is a partition $\gamma \in \sP$ that is an $S$-splitting of $\alpha$ and an $S$-splitting of $\beta$. Let $F, K$ be finite $S$-connected sets containing the identity such that $\alpha \le \beta^F$ and $\beta \le \alpha^K$. Thus $\alpha^K \le (\beta^{F})^K = \beta^{KF}$. Since $\beta$ is a coarsening of $\alpha^K$ and $KF$ is $S$-connected and contains the identity, the previous lemma implies $\gamma= \beta^{KF}$ is a splitting of $\alpha^K$, and therefore, is a splitting of $\alpha$. Of course, $\beta^{KF}$ is also a splitting of $\beta$.
\end{proof}

The next theorem explains how we will use this general theory to prove that $h(\Sigma,\alpha)$ does not depend on the choice of generating partition $\alpha$. 
\begin{thm}\label{thm:general}
Let $G$ be a countable group and suppose $S \subset G$ is a generating set for $G$. Let $(G,X,\mu)$ be a $G$-system. As above, let $\sP$ be the space of partitions of $X$ with finite entropy. Let $f:\sP \to \R\cup \{-\infty\}$ be an upper semi-continuous function such that $f$ is invariant under $S$-splittings (i.e., if $\alpha$ is a $S$-splitting of $\beta \in \sP$ then $f(\alpha) = f(\beta)$). Then for all generating partitions $\alpha, \beta \in \sP$, $f(\alpha)=f(\beta)$. 
\end{thm}

\begin{proof}
 Let $\alpha,\beta\in \sP$ be any two generating partitions. By theorem \ref{thm:dense}, there exists a sequence $\{\alpha'_i\}$ of partitions converging to $\beta$ such that each $\alpha'_i$ is $S$-equivalent to $\alpha$. By the previous proposition, $\alpha'_i$ and $\alpha$ have a common $S$-splitting. So $f(\alpha'_i)=f(\alpha)$ for all $i$. Since $f$ is upper semi-continuous, 
$$f(\alpha) = \lim_{i\to\infty} f(\alpha'_i) \le f(\beta).$$
The opposite inequality, $f(\beta)\le f(\alpha)$, is similar. This proves the theorem. 
\end{proof}

In the next section, we prove the inequality $H(\Sigma,\beta:F) - H(\beta) \ge H(\Sigma,\alpha:F) - H(\alpha)$ when $\alpha \ge \beta$. This will be useful in succeeding proofs. In section \ref{sec:uppersemi} we show that $H(\Sigma,\alpha:F)$ is well-defined when $\alpha$ is infinite and that $h(\Sigma,\cdot)$ is upper semi-continuous. In section \ref{sec:monotone}, we show that $h(\Sigma,\cdot)$ is invariant under $S$-splittings. Theorem \ref{thm:general} then implies theorem \ref{thm:K}. 

\section{A lower bound for the entropy of a coarsening}\label{sec:factor}

The proposition succeeding the two lemmas below will be used frequently in this paper.

\begin{lem}\label{lem:factor}
Let $(G,X,\mu)$ be a system with $G$ a countable group.  Let $\Sigma=\{\sigma_i\}$ be a map sequence for $G$ (i.e., a sequence of maps $\sigma_i:G \to \sym(m_i)$ with $m_i\to \infty$ as $i\to \infty$). Let $\alpha, \beta$ be finite partitions of $X$ and suppose that $\alpha$ refines $\beta$. Then for all finite $F\subset G$,
$$H(\Sigma,\beta:F)-H(\beta) \ge H(\Sigma,\alpha:F)-H(\alpha).$$
\end{lem}

\begin{proof}
Let $\alpha=(A_1,A_2,\dots, A_u)$ and $\beta=(B_1,B_2,\dots, B_v)$. Let $\epsilon >0$. Let $\sigma: G \to \sym(m)$. We will obtain an upper bound on the cardinality of $\FA(\sigma,\alpha:F,\epsilon)$ in terms of $|\FA\big(\sigma,\beta:F,\epsilon\big)|$.  

Let $b:\N \to \N$ be the map $b(i)=j$ if $A_i \subset B_j$. Define the coarsening map $\Phi:\FA(\sigma,\alpha:F,\epsilon) \to \FA(\sigma,\beta:F,\epsilon)$ as follows. If $\balpha=(\bA_1,\dots,\bA_u)$ then let $\Phi(\balpha)=\bbeta=(\bB_1,\dots,\bB_v)$ where $\bB_j = \bigcup_{i:b(i)=j} \bA_i$. We need to check that $d_F(\beta,\bbeta)\le \epsilon$ as claimed. 

Let $\zeta$ be the uniform probability measure on $\{1,\dots,m\}$. As in definition \ref{defn:FA}, given a function $\phi:F \to \N$, let $B_\phi = \bigcap_{f \in F} fB_{\phi(f)}$ and $\bB_\phi = \bigcap_{f\in F} \sigma(f)\bB_{\phi(f)}$. Define $A_\phi, \bA_\phi$ similarly.

Note that $B_\phi = \cup_\psi A_\psi$ where the union is over all $\psi:F \to \N$ such that $A_\psi \subset B_\phi$. Similarly, $\bB_\phi = \cup_\psi \bA_\psi$. Therefore,
\begin{eqnarray*}
d_F(\beta, \bbeta) &=& \sum_{\phi:F \to \N} \big|\mu(B_\phi) - \zeta(\bB_\phi)\big| \le \sum_{\psi:F \to \N} \big|\mu(A_\psi) - \zeta(\bA_\psi)\big| =d_F(\alpha, \balpha) \le \epsilon.
\end{eqnarray*}
This shows that $\bbeta \in \FA(\sigma,\beta:F,\epsilon)$ as claimed. 

Let $\bbeta=(\bB_1,\dots,\bB_v) \in \FA(\sigma,\beta:F,\epsilon)$. We will bound the cardinality of the inverse image $\Phi^{-1}(\bbeta)$. 

For a vector $\v=(v_1,v_2,\dots,v_u) \in \N^u$, let $\cA(\bbeta:\v)$ be the set of all $\balpha=(\bA_1,\dots,\bA_u) \in \Phi^{-1}(\bbeta)$ such that $|\bA_i|=v_i$ for all $1\le i \le u$. A partition $\balpha \in \cA(\bbeta:\v)$ is obtained from $\bbeta$ by subdividing each partition element $\bB_j$ into sets of cardinality $v_i$ for $b(i)=j$. Therefore,
$$|\cA(\bbeta:\v)| \le \Big(\prod_{j=1}^{v} |\bB_j|!\Big) \Big(\prod_{i=1}^u v_i!\Big)^{-1}.$$
If $\balpha \in \cA(\bbeta:\v)$ then, 
\begin{eqnarray}\label{eqn:bongo}
\sum_{i=1}^{u} |\mu(A_i) - \zeta(\bA_i)|=\sum_{i=1}^{u} |\mu(A_i)-v_i/m| \le d_F(\alpha,\balpha) \le \epsilon.
\end{eqnarray}
 So $v_i \approx \mu(A_i)m$. Stirling's approximation implies that there is a function $\delta:\R \to \R$ (depending only on $\alpha$ and $\beta$) such that if $m$ is sufficiently large then
$$|\cA(\bbeta:\v)| \le \exp\Big( \big(H(\alpha)-H(\beta) + \delta(\epsilon) \big)m \Big)$$ 
and $\lim_{\epsilon \to 0} \delta(\epsilon) = 0$.

By equation \ref{eqn:bongo}, the number of vectors $\v \in \N^u$ such that $\cA(\bbeta:\v)$ is nonempty is at most $(3m\epsilon)^u$. Thus,
\begin{eqnarray*}
|\Phi^{-1}(\bbeta)| &=& \sum_{\{\v | \cA(\bbeta:\v)\ne \emptyset\}} |\cA(\bbeta:\v)|\\ 
&\le & (3m\epsilon)^u\exp\Big( \big(H(\alpha)-H(\beta) + \delta(\epsilon) \big)m \Big) .
\end{eqnarray*}
Since $\bbeta \in \FA(\sigma,\beta:F,\epsilon)$ is arbitrary, this implies that if $m$ is sufficiently large,
\begin{eqnarray*}
\Big|\FA(\sigma,\alpha:F,\epsilon)\Big|
&\le & (3m\epsilon)^u \exp\Big( \big(H(\alpha)-H(\beta) + \delta(\epsilon) \big)m\Big) \Big|\FA(\sigma,\beta:F,\epsilon)\Big|.
\end{eqnarray*}
Since $\sigma$ is arbitrary, the definition of $H(\Sigma,\alpha:F,\epsilon)$ yields:
$$H\big(\Sigma,\alpha:F,\epsilon\big) \le  H(\alpha)-H(\beta) + \delta(\epsilon) +H\big(\Sigma,\beta:F,\epsilon\big).$$
Now take the limit as $\epsilon \to 0$ to obtain:
$$H(\Sigma,\alpha:F) \le H(\alpha)-H(\beta) + H(\Sigma,\beta:F).$$
\end{proof}

\begin{lem}\label{lem:welldefined}
Let $\alpha \in \sP$. Let $\{\alpha_n\}$ be a chain of $\alpha$ (as in definition \ref{defn:chain}). Let $\Sigma$ be a map sequence of $G$. Let $F\subset G$ a finite set. Then $\lim_{n\to\infty} H(\Sigma, \alpha_n:F)$ exists. Moreover, if $\{\beta_n\}$ is another chain of $\alpha$ with $\beta_n \ge \alpha_n$ for all $n$ then 
 $$\lim_{n\to\infty} H(\Sigma, \beta_n:F)  \le \lim_{n\to\infty} H(\Sigma, \alpha_n:F).$$
\end{lem}

\begin{proof}
By lemma \ref{lem:relative} and the previous lemma, for each $m > n \ge 1$,
$$H(\Sigma,\alpha_m:F) \le H(\Sigma,\alpha_n:F) + H(\alpha_m|\alpha_n).$$
Take the limsup as $m \to \infty$ to obtain
$$\limsup_{m \to \infty} H(\Sigma,\alpha_m:F) \le H(\Sigma,\alpha_n:F) + H(\alpha|\alpha_n)$$
for every $n \ge 1$. Take the liminf as $n \to \infty$ to obtain
$$\limsup_{m \to\infty} H(\Sigma,\alpha_m:F) \le \liminf_{n \to\infty} H(\Sigma,\alpha_n:F).$$
Here we are using $H(\alpha)<\infty$ which implies $H(\alpha|\alpha_n)$ tends to zero as $n\to \infty$. Hence the limit, $\lim_{n\to\infty} H(\Sigma,\alpha_n:F)$, exists.

If $\{\beta_n\}$ is another chain of $\alpha$ with $\beta_n \ge \alpha_n$ for all $n$ then 
$$H(\Sigma,  \beta_n:F) \le H(\Sigma, \alpha_n:F) + H(\beta_n|\alpha_n).$$
Since $H(\beta_n|\alpha_n) \le H(\alpha |\alpha_n) \to 0$ as $n\to \infty$, this implies that 
$$\lim_{n\to\infty} H(\Sigma,  \beta_n:F) \le \lim_{n\to\infty} H(\Sigma, \alpha_n:F).$$
\end{proof}

We can now remove the finiteness hypothesis in lemma \ref{lem:factor}.

\begin{prop}\label{prop:factor}
Let $G$ be a countable group, $\Sigma$ a map sequence of $G$ and $(G,X,\mu)$ a $G$-system. Let $\alpha, \beta \in \sP$ and suppose that $\alpha$ refines $\beta$. Then for all finite $F \subset G$,
$$H(\Sigma,\beta:F)-H(\beta) \ge H(\Sigma,\alpha:F)-H(\alpha).$$
\end{prop}

\begin{proof}
Let $\{\alpha_n\}$, $\{\beta_n\}$ be chains of $\alpha$ and $\beta$ respectively.  It follows from lemma \ref{lem:factor} that
$$H(\Sigma,\beta_n:F)-H(\beta_n) \ge H(\Sigma,\alpha_n \vee \beta_n:F)-H(\alpha_n \vee \beta_n)$$
for all $n$. Take the limit as $n\to\infty$ to obtain
\begin{eqnarray*}
\Big(\lim_{n\to\infty} H(\Sigma,\beta_n:F) \Big) -H(\beta) &\ge& \Big(\lim_{n\to\infty} H(\Sigma,\alpha_n \vee \beta_n:F)\Big)-H(\alpha)\\
&\ge& H(\Sigma,\alpha :F) -H(\alpha).
\end{eqnarray*}
Since this is true for every chain $\{\beta_n\}$ of $\beta$, the result follows.
\end{proof}

\section{Upper semicontinuity}\label{sec:uppersemi}


\begin{prop}\label{prop:semicontinuity}
Let $(G,X,\mu)$ be a system. Let $\Sigma=\{\sigma_i\}_{i=1}^\infty$ be a sequence of maps $\sigma_i:G \to \sym(m_i)$ with $m_i \to \infty$ as $i\to\infty$. Then for any finite $F \subset G$ the map $H(\Sigma,\cdot:F):\sP \to \R \cup \{-\infty\}$ is upper semi-continuous.
\end{prop}

\begin{proof}
Let $\alpha \in \sP$ and let $\{\beta^i\}_{i=1}^\infty \subset \sP$ be a sequence converging to $\alpha$. It suffices to show that $\limsup_i H(\Sigma,\beta^i:F) \le H(\Sigma,\alpha:F)$. Because $H(\Sigma,\beta^i:F)$ does not depend on the ordering of the atoms of $\beta^i$, we may assume, without loss of generality, that if $\alpha=(A_1, A_2, \dots)$ and $\beta^i=(B^i_1, B^i_2, \dots)$ then 
$$\lim_{i\to\infty} \mu( A_j \Delta B^i_{j}) = 0$$
for all $j$.

Let us now assume that $\beta^i=(B^i_1,\dots,B^i_u)$ and $\alpha=(A_1,\dots,A_u)$ are finite partitions with the same number of atoms. As in definition \ref{defn:FA}, given a function $\phi:F \to \N$, let $A_\phi = \bigcap_{f \in F} fA_{\phi(f)}$ and $B^i_\phi = \bigcap_{f\in F} fB^i_{\phi(f)}$. Define
$$d_F(\alpha,\beta^i)=\sum_{\phi:F \to \N} \Big| \mu( A_\phi) - \mu(B^i_\phi )\Big|.$$
Then $d_F(\alpha,\beta^i)\to 0$ as $i\to\infty$ since $\beta^i \to \alpha$ (see lemma \ref{lem:alpha^n}). An elementary computation shows for any map $\sigma:G \to \sym(m)$, any $\epsilon>0$ and any $i\ge 0$,
$$\FA\big(\sigma, \alpha:F,\epsilon + d_F(\alpha,\beta^i)\big) \supset \FA(\sigma, \beta^i:F,\epsilon).$$
Therefore if $c>1$ is arbitrary,
$$H(\Sigma,\alpha:F, c d_F(\alpha,\beta^i) ) \ge H(\Sigma,\beta^i:F).$$
Take the limsup as $i$ tends to infinity to obtain
$$H(\Sigma,\alpha:F ) \ge \limsup_{i\to\infty} H(\Sigma,\beta^i:F)$$
as claimed.

Now let $\beta^i=(B^i_1,B^i_2,\dots)$ and $\alpha=(A_1,A_2,\dots)$ be infinite partitions. It is allowed that some of the atoms of $\alpha$ and/or $\beta^i$ are empty. So this case includes the finite case. Let $\{\alpha_n\}_{n=1}^\infty$ be a chain of $\alpha$. For each $n$, there exists a finite partition $\pi_n=\{P_{1,n}, \dots, P_{c_n,n}\}$ of $\N$ such that every atom of $\alpha_n$ is of the form $\bigcup_{j \in P_{k,n}} A_j$ for some $k$. Let $\beta^i_{n}$ be the corresponding coarsening of $\beta^i$. That is, let $\beta^i_{n}$ be the partition whose atoms are all of the form $\bigcup_{j \in P_{k,n}} B^i_{j}$ (for $1\le k \le c_n$).


Proposition \ref{prop:factor} implies that
$$H(\Sigma, \beta^i:F) - H(\beta^i| \beta_n^i) \le H(\Sigma, \beta^i_n:F)$$
for all $i,n$. The previous case implies
$$\limsup_i H(\Sigma,\beta^i_n:F) \le H(\Sigma,\alpha_n:F)$$
for all $n$. Thus,
$$\limsup_i H(\Sigma,\beta^i:F) \le H(\Sigma,\alpha_n:F) + \limsup_i H(\beta^i| \beta_n^i) = H(\Sigma, \alpha_n:F) + H(\alpha|\alpha_n).$$
Since this is true for all $n$ and for every chain $\{\alpha_n\}$ of $\alpha$, it follows that
$$\limsup_{i \to\infty}  H(\Sigma,\beta^i:F) \le H(\Sigma,\alpha:F)$$
as claimed.
\end{proof}

\begin{prop}\label{prop:welldefined}
Let $\alpha \in \sP$ and $\{\alpha_n\}$ a chain of $\alpha$. Let $\Sigma$ be a map sequence of $G$. Let $F\subset G$ a finite set. Then 
$$ H(\Sigma,\alpha:F)=\lim_{n\to\infty} H(\Sigma, \alpha_n:F).$$
\end{prop}

\begin{proof}
It follows from the definition that $H(\Sigma,\alpha:F) \le \lim_{n\to\infty} H(\Sigma,\alpha_n:F)$. Because $\alpha_n \to \alpha$ as $n\to\infty$, the previous proposition implies that $H(\Sigma,\alpha:F) \ge \lim_{n\to\infty} H(\Sigma,\alpha_n:F)$.
\end{proof}

\begin{cor}\label{cor:uppersemicontinuous}
Let $(G,X,\mu)$ be a system. Let $\Sigma$ be a map sequence of $G$. Then the map $h(\Sigma,\cdot):\sP \to \R \cup \{-\infty\}$ is upper semi-continuous.
\end{cor}

\begin{proof}
By definition $h(\Sigma,\cdot) = \inf_F H(\Sigma,\cdot:F)$ where the infimum is over all finite $F \subset G$. Since an infimum of upper semi-continuous functions is upper semi-continuous, proposition \ref{prop:semicontinuity} implies this corollary.
\end{proof}

\section{Monotonicity}\label{sec:monotone}

In this section, we prove that $H(\Sigma,\alpha:F)$ is monotone decreasing under $S$-splittings. We handle the finite partition case first.

\begin{prop}\label{prop:monotone}
Let $F \subset G$ be finite and suppose $e \in F$. If $\alpha$ is an $F$-splitting of a finite partition $\beta \in \sP$ and $\Sigma=\{\sigma_i\}_{i=1}^\infty$ is a sequence of maps $\sigma_i:G \to \sym(m_i)$ with $m_i \to \infty$ then $H(\Sigma, \alpha:F)\le H(\Sigma,\beta:F)$.
\end{prop}

\begin{proof}
Intuitively, the proposition is true because any approximation to $\beta$ on $\{1,\dots,m\}$ can be split into an approximation of $\alpha$ and approximately all of the approximations to $\alpha$ are 
obtained this way. The proof is a matter of making this intuition precise.

Let $1/4>\epsilon>0$ and $\sigma: G \to \sym(m)$ . We will obtain an upper bound on the cardinality of $\FA(\sigma,\alpha:F,\epsilon)$ in terms of $|\FA\big(\sigma,\beta:F,\epsilon\big)|$.

It suffices to consider only the special case in which $\alpha=(A_1,\dots,A_u)$ is a simple $F$-splitting of $\beta=(B_1,\dots,B_v)$. So there exists $f\in F$ and a coarsening $\xi=(X_1,\dots,X_w)$ of $\beta$ such that $\alpha = \beta \vee f\xi$.  

Let $b,x:\{1,\dots,u\} \to \N$ be the maps defined by $$A_i = B_{b(i)} \cap fX_{x(i)}.$$

Define the ``coarsening'' map $\Phi:\FA(\sigma,\alpha:F,\epsilon) \to \FA\big(\sigma,\beta:F,\epsilon\big)$ as follows. For $\balpha=(\bA_1, \dots, \bA_u) \in \FA(\sigma,\alpha:F,\epsilon)$, let $\Phi(\balpha)=\bbeta=(\bB_1,\dots,\bB_v)$ where $\bB_j=\bigcup_{i:b(i)=j} \bA_i$.
As in the proof of lemma \ref{lem:factor}, $d_F(\bbeta,\beta)\le \epsilon$ so $\bbeta \in \FA(\sigma,\beta:F,\epsilon)$ as claimed.

Next, we obtain an upper bound on the cardinality of $\Phi^{-1}(\bbeta)$ where $\bbeta \in \FA\big(\sigma,\beta:F,\epsilon\big)$ is a fixed partition. This bound will not depend on the choice of $\bbeta$ and thus, we will be able to use it to bound $|\FA(\sigma,\alpha:F,\epsilon)|$. In order to obtain the bound, we will show that every partition in $\Phi^{-1}(\bbeta)$ is ``close'' to a particular partition which we construct directly from $\bbeta$ next.

Let $t:\{1,\dots,v\} \to \N$ be defined by $B_i \subset X_{t(i)}$. Define $\bxi =(\bX_1,\dots,\bX_w)$ by 
$$\bX_j = \bigcup_{i:t(i)=j} \bB_i.$$
Let $\balpha=(\bA_1,\dots,\bA_u)$ be defined by 
$$\bA_i = \bB_{b(i)} \cap \sigma(f) \bX_{x(i)}.$$

Let $\zeta$ be the uniform probability measure on $\{1,\dots,m\}$. Everything we need to know about $\balpha$ is contained in the claims below.

{\bf Claim 1}: if $\delta=(D_1,\dots,D_u) \in \Phi^{-1}(\bbeta)$, then 
$$\zeta\Big( \bigcup_{i=1}^u D_i \Delta \bA_i\Big) \le \epsilon.$$

{\bf Claim 2}: if $\delta$ is as above and $z \in D_i - \bA_i$ for some $i$ then $z \in D_i \cap \sigma(f)D_j$ for some $1 \le j \le u$ such that $A_i \cap fA_j = \emptyset$.

Let us see how claim 2 implies claim 1. First note that $\bigcup_{i=1}^u D_i \Delta \bA_i = \bigcup_{i=1}^{u} D_i - \bA_i$. To see this, let $z \in \bigcup_{i=1}^u D_i \Delta \bA_i$. Then there is an $j$ such that $z \in D_j \Delta \bA_j$. If $z \in D_j -\bA_j$ then $z\in \bigcup_{i=1}^{u} D_i - \bA_i$. Otherwise $z\in \bA_j - D_j$. Since $\delta$ is a partition, there is a $k\ne j$ such that $z \in D_k$. But then $z \in D_k - \bA_k$ since $\bA_k$ does not intersect $\bA_j$. Since $z$ is arbitrary, this shows that $\bigcup_{i=1}^u D_i \Delta \bA_i \subset \bigcup_{i=1}^{u} D_i - \bA_i$. The opposite inclusion is trivial.

Claim 2 implies that
\begin{eqnarray*}
\zeta\Big(\bigcup_{i=1}^u D_i \Delta \bA_i\Big)= \zeta\Big(  \bigcup_{i=1}^{u} D_i - \bA_i\Big)\le  \sum_{i,j=1}^u \Big|\mu \big(A_i \cap fA_j\big) - \zeta\big(D_i \cap \sigma(f)D_j\big)\Big| &\le& d_F(\alpha,\delta)\le \epsilon.
\end{eqnarray*}
This proves claim 1. To prove claim 2, let $z \in D_i- \bA_i$ for some $1\le i \le u$. By definition, $\bA_i = \bB_{b(i)} \cap \sigma(f)\bX_{x(i)}$. Since $\Phi(\delta)=\bbeta$, $z \in D_i \subset \bB_{b(i)}$. So, it must be that $z \notin \sigma(f)\bX_{x(i)}$. Let $j$ be such that $z \in \sigma(f)D_j \subset \sigma(f)\bB_{b(j)} \subset \sigma(f)\bX_{t(b(j))}$. Note $t(b(j)) \ne x(i)$.

By definition,
$$A_i \cap fA_j \subset  fX_{x(i)} \cap fX_{t(b(j))}.$$ 
Since $t(b(j))  \ne x(i)$, $X_{t(b(j))} \cap X_{x(i)}$ is empty. So $A_i \cap fA_j$ is empty. This proves claim 2.

It follows from claim 1 that any partition in $\Phi^{-1}(\bbeta)$ can be obtained from $\balpha$ by ``relabeling'' a subset of $\{1,\dots,m\}$ of cardinality at most $\epsilon m$. That is, if $\delta \in \Phi^{-1}(\bbeta)$, then there exists a set of cardinality $n  = \lceil \epsilon m \rceil$ in $\{1,\dots,m\}$ such that $\delta$ can be obtained from $\balpha$ by redefining $\balpha$ on this set. Since $\epsilon < 1/4$ this implies
$$|\Phi^{-1}(\bbeta)| \le {m \choose n}u^{n}.$$
If $m$ is sufficiently large then Stirling's approximation implies
$$|\Phi^{-1}(\bbeta)| \le \exp\big(H(2\epsilon,1-2\epsilon)m + 2\epsilon m\log(u)\big)$$
where $H(x,y)=-x\log(x)-y\log(y)$. Since $\bbeta \in \FA\big(\sigma,\beta:F,\epsilon\big)$ is arbitrary, this implies that if $m$ is sufficiently large, 
$$\big|\FA(\sigma,\alpha:F,\epsilon)\big| \le \exp\big(H(2\epsilon,1-2\epsilon)m + 2\epsilon m\log(u)\big)\big|\FA\big(\sigma,\beta:F,\epsilon\big)\big|.$$
Thus,
$$H\big(\Sigma,\alpha:F,\epsilon\big) \le H(2\epsilon,1-2\epsilon) + 2\epsilon \log(u)+H\big(\Sigma,\beta:F,\epsilon\big).$$
Now take the limit as $\epsilon \to 0$ to obtain $H(\Sigma,\alpha:F) \le H(\Sigma,\beta:F)$, as claimed.
\end{proof}

Next, we remove the finiteness assumption in the above proposition.
\begin{prop}\label{prop:monotone2}
Let $F \subset G$ be finite and suppose $e \in F$. If $\alpha$ is an $F$-splitting of a partition $\beta \in \sP$ and $\Sigma=\{\sigma_i\}_{i=1}^\infty$ is a map sequence then $H(\Sigma, \alpha:F)\le H(\Sigma,\beta:F)$.
\end{prop}

\begin{proof}
It suffices to assume $\alpha$ is a simple $F$-splitting of $\beta$. So there is an $ f\in F$ and a coarsening $\xi$ of $\beta$ such that $\alpha = \beta \vee f\xi$. 

Let $\{\beta_n\}$, $\{\xi_n\}$ be chains for $\beta$ and $\xi$ respectively such that $\xi_n$ coarsens $\beta_n$ for all $n$. Then $\{\beta_n \vee f\xi_n\}$ is a chain for $\alpha$. The previous proposition implies
$$H(\Sigma, \beta_n \vee f\xi_n:F) \le H(\Sigma, \beta_n:F)$$
for all $n$. Proposition \ref{prop:welldefined} now implies that $H(\Sigma,\beta \vee f\xi:F) \le H(\Sigma,\beta:F)$ as claimed.
\end{proof}

\begin{prop}\label{prop:monotone3}
Let $S \subset G$. Let $\Sigma$ be a sofic approximation to $G$. If $\alpha$ is an $S$-splitting of $\beta \in \sP$ then $h(\Sigma,\alpha)=h(\Sigma,\beta)$.
\end{prop}

\begin{proof}
 We may assume that $\alpha$ is a simple $S$-splitting of $\beta$. So there exists an $t \in S$ and a coarsening $\xi$ of $\beta$ such that $\alpha = \beta \vee t\xi$. Note that $\beta \vee t\beta$ is an $S$-splitting of $\alpha$ by lemma \ref{lem:$S$-splittings}. So the previous proposition implies 
$$h(\Sigma, \beta \vee t\beta) \le h(\Sigma,\alpha) \le h(\Sigma,\beta).$$
So it suffices to show that $h(\Sigma,\beta \vee t\beta) \ge h(\Sigma,\beta)$. 

Let us assume for now that $\beta=(B_1,\dots,B_v)$ is finite. Let $\alpha:=\beta \vee t\beta=(A_1, A_2, \dots, A_u)$. Let $x,y:\{1,\dots,u\} \to \N$ be maps determined by $A_i = B_{x(i)} \cap tB_{y(i)}$. We assume that the map $i \mapsto (x(i),y(i))$ surjects onto $\{1,\ldots,v\} \times \{1,\ldots,v\}$. Thus some of the $A_i$'s may be empty.

Let $F \subset G$ be a finite set with $e,t \in F$ and let $\epsilon,\delta>0$. Let $\sigma:G \to \sym(m)$ be an $(F, \delta)$-approximation to $G$. Let $V(F) \subset \{1,\dots,m\}$ be the set of all elements $v$ such that for all $f_1,f_2 \in F$, 
$$\sigma(f_1)\sigma(f_2)v = \sigma(f_1f_2)v$$
and $\sigma(f_1)v \ne \sigma(f_2)v$ if $f_1 \ne f_2$. Since $\sigma$ is an $(F,\delta)$-approximation, $|V(F)| \ge (1-\delta)m$.

Define 
$$\Psi: \FA\big(\sigma, \beta: F \cup Ft, \epsilon\big) \to \FA\big(\sigma, \beta \vee t\beta:F, \epsilon+5|F|\delta\big)$$
by $\Psi(\bbeta) = \balpha = (\bA_1,\dots,\bA_u)$ where if $\bbeta=(\bB_1,\dots,\bB_v)$ then $\bA_i = \bB_{x(i)} \cap \sigma(t)\bB_{y(i)}$. Note that $\Psi$ is injective.

It is implicitly claimed above that $d_F(\alpha,\Psi(\bbeta))=d_F(\alpha,\balpha)\le \epsilon + 5|F|\delta$. Let us check this. As usual, for $\phi:F \to \N$ set $A_\phi = \bigcap_{f \in F} fA_{\phi(f)}$. Similar formulas apply to $\bA_\phi$ and to $B_\psi, \bB_\psi$ for $\psi:F\cup Ft \to \N$. 

Let $\cG$ be the set of all functions $\phi: F \to \N$ such that for every $f\in F$ with $ft \in F$, $x(\phi(ft))=y(\phi(f))$. Observe that if $\phi\notin \cG$ and $f\in F$ is such that $ft \in F$ but $x(\phi(ft))\ne y(\phi(f))$ then
$$A_\phi \subset ftA_{\phi(ft)} \cap fA_{\phi(f)} \subset ftB_{x(\phi(ft))} \cap ftB_{y(\phi(f))} = \emptyset.$$
Similarly,
$$\bA_\phi \subset \sigma(ft)\bB_{x(\phi(ft))} \cap \sigma(f)\sigma(t)\bB_{y(\phi(f))} \subset \sigma(f)\sigma(t)V(F)^c$$
where $V(F)^c$ denotes the complement of $V(F)$. Since $\zeta(V(F)^c) \le \delta$ and $\{\bA_\phi\}$ is pairwise disjoint,
\begin{eqnarray*}
d_F(\alpha, \balpha) = \sum_{\phi:F \to \N} \Big| \mu( A_\phi ) - \zeta( \bA_{\phi} )\Big|\le |F|\delta+ \sum_{\phi \in \cG}  \Big| \mu( A_\phi ) - \zeta( \bA_{\phi} )\Big|.
\end{eqnarray*}

On the other hand, if $\phi\in \cG$ then define $\tphi:F \cup Ft \to \N$ as follows. For $f\in F$, let $\tphi(f)=x(\phi(f))$ and $\tphi(ft)=y(\phi(f))$. This is well-defined by the definition of $\cG$. Observe that
$$A_\phi = \bigcap_{f \in F} fA_{\phi(f)} = \bigcap_{f \in F} fB_{x(\phi(f))} \cap ftB_{y(\phi(f))} =\bigcap_{f \in F \cup Ft} fB_{\tphi(f)}=B_\tphi.$$

It is almost true that $\bA_\phi = \bB_\tphi$. If $\sigma$ is not a homomorphism then this equality can fail. To be precise, observe that
$$\bA_\phi = \bigcap_{f\in F} \sigma(f)\bA_{\phi(f)} = \bigcap_{f\in F} \sigma(f)\bB_{x(\phi(f))} \cap \sigma(f)\sigma(t) \bB_{y(\phi(f))},$$
while
$$\bB_\tphi = \bigcap_{f\in F} \sigma(f)\bB_{x(\phi(f))} \cap \sigma(ft) \bB_{y(\phi(f))}.$$
Hence 
$$\bA_\phi \Delta \bB_\tphi \subset \bigcup_{f\in F} \sigma(ft)V(F)^c \cup \sigma(f)\sigma(t)V(F)^c$$
where $V(F)^c$ denotes the complement of $V(F)$. 

Each of the collections $\{\bA_\phi - \bB_\tphi\}_{\phi \in \cG}$ and $\{\bB_\tphi - \bA_\phi\}_{\phi \in \cG}$ is pairwise disjoint. Thus,
$$\sum_{\phi \in \cG} \Big|\bA_\phi \Delta \bB_\tphi\Big| \le 2\Big|\bigcup_{f\in F} \sigma(ft)V(F)^c \cup \sigma(f)\sigma(t)V(F)^c\Big| \le 4|F|\delta m.$$
This implies $\sum_{\phi \in \cG}|\zeta( \bA_\phi ) - \zeta( \bB_{\tphi} )|\le 4|F|\delta$. So,
\begin{eqnarray*}
d_F(\alpha, \balpha) &=& \sum_{\phi:F \to \N} \Big| \mu( A_\phi ) - \zeta( \bA_{\phi} )\Big|\\
&\le& |F|\delta+ \sum_{\phi \in \cG}\Big| \mu( A_\phi ) - \zeta( \bA_{\phi} )\Big|\\
&\le&|F|\delta + \sum_{\phi \in \cG} \Big| \mu(B_{\tphi}) - \zeta( \bB_{\tphi})\Big| + \Big| \zeta( \bA_\phi ) - \zeta( \bB_{\tphi} )\Big|\\
&\le& d_{F\cup Ft}(\beta,\bbeta) + 5|F|\delta \le \epsilon + 5|F|\delta.
\end{eqnarray*}
This proves the claim. 

Since $\Psi$ is injective,
$$\Big|\FA\big(\sigma, \beta: F \cup Ft, \epsilon\big)\Big| \le \Big|\FA\big(\sigma, \beta \vee t\beta:F,\epsilon + 5|F|\delta\big)\Big|.$$
Since $\Sigma=\{\sigma_i\}$ is a sofic approximation, each $\sigma_i$ is a $(F,\delta_i)$-approximation to $G$ for some $\delta_i\ge 0$ with $\delta_i \to 0$ as $i\to\infty$. Thus  if $c>1$ is arbitrary then
$$H\big(\Sigma, \beta: F \cup Ft, \epsilon\big) \le H\big(\Sigma, \beta \vee t\beta:F,c\epsilon\big).$$
Let $\epsilon \to 0$ to obtain
\begin{eqnarray}\label{eqn:thingaling}
H\big(\Sigma, \beta: F \cup Ft\big) \le H\big(\Sigma, \beta \vee t\beta:F\big).
\end{eqnarray}
Finally, take the infimum over all finite sets $F\subset G$ with $e,t \in F$ to obtain
$$h(\Sigma, \beta) \le h(\Sigma,\beta \vee t\beta).$$
This finishes the proof in the case that $\beta$ is finite.

Now suppose that $\beta=(B_1,\dots) \in \sP$ is a possibly infinite partition. Let $\{\beta_n\}$ be a chain of $\beta$. For every finite $F \subset G$ with $e,t \in F$, equation \ref{eqn:thingaling} implies
$$H(\Sigma, \beta_n :F\cup Ft) \le H(\Sigma, \beta_n \vee t\beta_n:F).$$
Since $\{\beta_n \vee t\beta_n\}$ is a chain for $\beta \vee t\beta$, proposition \ref{prop:welldefined} implies
$$H(\Sigma, \beta :F\cup Ft) \le H(\Sigma, \beta \vee t\beta:F).$$
Now take the infimum over finite sets $F \subset G$ with $e,t \in F$ to obtain
$$h(\Sigma, \beta) \le h(\Sigma,\beta \vee t\beta).$$
This finishes the proof.
\end{proof}

We can now prove theorem \ref{thm:K}.

\begin{proof}[Proof of theorem \ref{thm:K}]
By corollary \ref{cor:uppersemicontinuous}, $h(\Sigma,\cdot):\sP \to \R \cup \{-\infty\}$ is upper semi-continuous. By the previous proposition, $h(\Sigma,\cdot)$ is invariant under splittings. The theorem now follows from theorem \ref{thm:general}.
\end{proof}

\section{Bernoulli shifts over a sofic group}\label{sec:Bernoulli}
In this section, we calculate the entropy of a product system in which one of the factors is Bernoulli. To be precise we need the following.
\begin{defn}
Let $(G,X,\mu)$ and $(G,Y,\nu)$ be two systems. Then the {\bf product system} $(G, X \times Y, \mu \times \nu)$ is defined by $g(x,y)=(gx, gy)$. If $\alpha$ is a partition of $X$ and $\beta$ is a partition of $Y$ then $\alpha \times \beta :=\{A \times B~|~A \in \alpha, B \in \beta\}$ is a partition of $X \times Y$. 
\end{defn}

\begin{thm}\label{thm:product}
Let $G$ be a group with sofic approximation $\Sigma$. Let $K$ be a finite or countably infinite set and let $\kappa$ be a probability measure on $K$ such that $H(\kappa)<\infty$. Let $\beta$ be the canonical partition of $K^G$. I.e., $\beta=\{B_k~|~k\in K\}$ where $B_k=\{x \in K^G~|~x(e)=k\}$.

If $(G,X,\mu)$ is any $G$-system and $\alpha$ is a partition of $X$ with $H(\alpha)<\infty$, then for every finite $F \subset G$,
$$H\big(\Sigma, \alpha \times \beta:F \big) = H(\Sigma,\alpha:F) + H(\kappa).$$
This implies:
$$h\big(\Sigma, G,X\times K^G,\mu \times \kappa^G) = h(\Sigma, G,X,\mu) + H(\kappa)$$
if $(G,X,\mu)$ has a finite-entropy generating partition. In particular, $h(\Sigma,G,K^G,\kappa^G)=H(\kappa)$.

\end{thm}

\begin{proof}
We will first obtain the upper bound. Let $\tau=\{K^G\}$ be the trivial partition of $K^G$. Proposition \ref{prop:factor} implies
$$H(\Sigma,\alpha \times \tau:F) -H(\alpha\times \tau) \ge H(\Sigma, \alpha \times \beta:F) - H(\alpha \times \beta).$$
It is a standard exercise to prove that $H(\alpha \times \beta) = H(\alpha) + H(\beta)=H(\alpha)+H(\kappa)$. Since $H(\alpha \times\tau)=H(\alpha)$, 
$$H(\Sigma, \alpha \times \beta:F) \le H(\Sigma,\alpha:F) + H(\kappa).$$
This proves the upper bound.

Observe that, by taking coarsenings of $\alpha$ and $\beta$, proposition \ref{prop:welldefined} implies that we may assume, without loss of generality, that $\alpha$ and $\beta$ are finite. Of course, this means that $K$ is finite. So let $p$ be the number of atoms of $\alpha$ and let $q=|K|$ be the number of atoms of $\beta$.

To begin, let $\epsilon, \delta>0$, be such that $2\delta< (pq)^{-|F|}\epsilon$. Let $\sigma: G \to \sym(m)$ be an $(F,\delta)$-approximation to $G$. Let $V(F) \subset \{1,\dots,m\}$ be the set of all elements $v$ such that for all $f_1,f_2 \in F$, 
$$\sigma(f_1)\sigma(f_2)v = \sigma(f_1f_2)v$$
and $\sigma(f_1)v \ne \sigma(f_2)v$ if $f_1 \ne f_2$. By definition, $|V(F)| \ge (1-\delta)m$.

Let $K=\{1,\dots,q\}$, $B_i=\{x\in K^G~|~x(e)=i\}$, $\beta=(B_1,\dots,B_q)$. Identify $K^m$ with the set of all ordered partitions $\bbeta=(\bB_1,\dots,\bB_q)$ of $\{1,\dots,m\}$ via the map $(x_1,\dots,x_m) \in K^m \mapsto (\bB_1,\dots,\bB_q)$ where 
$$\bB_i = \big\{ j \in \{1,\dots,m\}~|~ x_j=i\big\}.$$

Let $\kappa^m$ be the product measure on $K^m$, $\bbeta=(\bB_1,\dots,\bB_q)$ be a random element of $K^m$ with law $\kappa^m$ and $\balpha=(\bA_1,\dots,\bA_p) \in \FA(\sigma,\alpha:F,\epsilon)$. We will estimate the probability that $d_F( \alpha \times \beta, \balpha \vee \bbeta)\le 2\epsilon.$ Of course, $d_F( \alpha \times \beta, \balpha \vee \bbeta)$ depends on a choice of ordering of the two partitions. But we claim that the partitions $\balpha \vee \bbeta$ and $\alpha \times \beta$ can be ordered canonically so that
\begin{eqnarray}\label{eqn:d_F}
d_F( \alpha \times \beta, \balpha \vee \bbeta)=\sum_{\phi,\psi:F \to \N} \Big| \mu\times \kappa^G(A_\phi \times B_\psi) - \zeta(\bA_\phi \cap \bB_\psi)\Big|.
\end{eqnarray}
A word about the notation used above is in order. As in definition \ref{defn:FA}, if $\phi:F \to \N$ then $A_\phi = \bigcap_{f \in F} fA_{\phi(f)}$, $\bA_\phi = \bigcap_{f\in F} \sigma(f)\bA_{\phi(f)}$ and similar formulas hold for $B_\psi, \bB_\psi$. Note that
$$\bigcap_{f \in F} f \big(A_{\phi(f)} \times B_{\psi(f)}\big) = A_\phi \times B_\psi,$$
$$\bigcap_{f\in F} f\big(\bA_{\phi(f)} \cap \bB_{\psi(f)} \big) = \bA_\phi \cap \bB_\psi.$$
This justifies equation \ref{eqn:d_F}.

In order to estimate the probability that $d_F(\balpha \vee \bbeta, \alpha \times \beta)\le 2\epsilon$, fix functions $\phi,\psi:F \to \N$. For each $v \in \{1,\dots,m\}$, let $Z_v$ equal $1$ if $v \in V(F) \cap \bA_\phi \cap \bB_\psi$ and $0$ otherwise. Let $Z=\sum_{v=1}^m Z_v = |V(F) \cap \bA_{\phi} \cap \bB_\psi|$.

If $v \notin V(F) \cap \bA_\phi$ then the expected value of $Z_v$, denoted $\bE[Z_v]$, equals $0$. Otherwise, since $\sigma(f_1)v \ne \sigma(f_2)v$ for $v\in V(F)$ if $f_1 \ne f_2 \in F$, it follows that $\bE[Z_v] = \kappa^G(B_\psi)$. Thus
\begin{eqnarray*}
\bE\big[  |V(F) \cap \bA_{\phi} \cap \bB_\psi| \big] = \sum_{v=1}^m \bE[Z_v]= \kappa^G(B_\psi)\big|V(F) \cap \bA_{\phi}\big| .
\end{eqnarray*}

To estimate the variance of $Z$, denoted $\Var(Z)$, we will first bound $\bE[Z_vZ_w]$ for $1\le v,w \le m$. If either $v$ or $w$ is not in  $V(F) \cap \bA_\phi$ then $\bE[Z_vZ_w]=Z_vZ_w=0$. If $\sigma(f_1)v \ne \sigma(f_2)w$ for any $f_1,f_2 \in F$ then $Z_v$ and $Z_w$ are independent, in which case $\bE[Z_vZ_w]=\bE[Z_v]\bE[Z_w]$. So the number of non-independent pairs $(Z_v,Z_w)$ with both $v,w \in V(F) \cap \bA_\phi$ is at most $| V(F) \cap \bA_\phi| |F|^2$. So
\begin{eqnarray*}
\Var(Z) &=& -\bE[Z]^2 + \bE[Z^2]\\
&=& -\bE[Z]^2 + \sum_{v,w=1}^m \bE[Z_vZ_w]\\
&\le&  -\bE[Z]^2 + | V(F) \cap \bA_\phi| |F|^2 + \bE[Z]^2 \le m|F|^2.
\end{eqnarray*}
Chebyshev's inequality applied to $Z/m$ implies that for any $t >0$,
\begin{eqnarray}\label{eqn:cheby}
Pr\Big[ \big|Z/m - \bE[Z]/m \big| > t \Big] &\le& \frac{\Var(Z)}{m^2t^2} \le \frac{|F|^2}{m t^2}.
\end{eqnarray}
Here $Pr[\cdot ]$ denotes ``the probability that''. Observe that 
$$Z/m = |V(F)\cap \bA_\phi \cap \bB_\psi|/m = \zeta\big(V(F) \cap \bA_\phi \cap \bB_\psi\big),$$
$$\bE[Z]/m = \kappa^G(B_\psi)|V(F)\cap \bA_\phi|/m = \zeta\big(V(F) \cap \bA_\phi\big)\kappa^G(B_\psi).$$
Since $|V(F)| \ge (1-\delta)m$, it follows that
$$\Big| \zeta\big(V(F) \cap \bA_\phi \cap \bB_\psi\big) - \zeta(\bA_\phi \cap \bB_\psi) \Big| \le \delta,$$
$$\Big|  \zeta\big(V(F) \cap \bA_\phi\big)\kappa^G(B_\psi) - \zeta(\bA_\phi)\kappa^G(B_\psi) \Big| \le \delta.$$
So equation \ref{eqn:cheby} implies
\begin{eqnarray*}\label{eqn:cheby2}
Pr\Big[ \big|\zeta(\bA_\phi \cap \bB_\psi) - \zeta(\bA_\phi)\kappa^G(B_\psi)  \big| > t \Big] &\le&  \frac{|F|^2}{m (t-2\delta)^2}
\end{eqnarray*}
for any $t> 2\delta$. Set $t = (pq)^{-|F|}\epsilon$ ($p$ is the number of atoms of $\alpha$ and $q=|K|$ is the number of atoms of $\beta$). By choice of $\delta$, $t > 2\delta$. Then if $m$ is sufficiently large,
\begin{eqnarray*}\label{eqn:cheby3}
Pr\Big[ \big|\zeta(\bA_\phi \cap \bB_\psi) - \zeta(\bA_\phi)\kappa^G(B_\psi)  \big| > (pq)^{-|F|}\epsilon \Big] &\le&  (pq)^{-|F|}\epsilon.
\end{eqnarray*}
Therefore, the probability that $|\zeta(\bA_\phi \cap \bB_\psi) - \zeta(\bA_\phi)\kappa^G(B_\psi)  | > (pq)^{-|F|}\epsilon$ for {\it some} $\phi$ and $\psi$ is at most $\epsilon$. If this event does not occur then
\begin{eqnarray*}\label{eqn:d_F2}
d_F(\alpha \times \beta, \balpha \vee \bbeta)&=&\sum_{\phi,\psi:F \to \N} \big| \mu\times \kappa^G(A_\phi \times B_\psi) - \zeta(\bA_\phi \cap \bB_\psi)\big|\\
&\le&\sum_{\phi,\psi:F \to \N} \big| \mu(A_\phi)\kappa^G(B_\psi) - \zeta(\bA_\phi)\kappa^G(B_\psi)\big| +  (pq)^{-|F|}\epsilon\\
&\le& \epsilon+ \sum_{\phi:F \to \N} \big| \mu(A_\phi) - \zeta(\bA_\phi)\big|\\
&=&\epsilon + d_F(\alpha,\balpha) \le 2\epsilon.
\end{eqnarray*}
So if $m$ is sufficiently large, then
$$Pr\Big[ d_F(\balpha \vee \bbeta, \alpha \times \beta) \le 2\epsilon \Big] \ge 1-\epsilon.$$

It follows from the Shannon-McMillan theorem (or, more simply, from the law of large numbers) that for all $\epsilon>0$ there exists an $M>0$ such that if $m>M$ then there is a set $Q \subset K^m$ such that 
\begin{itemize}
\item $\kappa^m(Q) > 1-\epsilon$,
\item for all $q \in Q$, $\exp\big(-H(\kappa)m - \epsilon m\big) \le \kappa^m\big(\{q\}\big) \le \exp\big(-H(\kappa)m + \epsilon m\big)$.
\end{itemize}

Let $Q_0$ be the set of all $\bbeta \in Q$ such that $d_F(\balpha \vee \bbeta, \alpha \times \beta) \le 2\epsilon$. Then $\kappa^m(Q_0) \ge 1- 2\epsilon$ and
$$|Q_0| \ge \kappa^m(Q_0)\exp(H(\kappa)m -\epsilon m) \ge (1-2\epsilon)\exp(H(\kappa)m - \epsilon m)$$
for all sufficiently large $m$. Since this is true for every $\balpha \in \FA(\sigma, \alpha:F,\epsilon)$ it follows that
$$\big|\FA(\sigma,\alpha\times\beta:F,2\epsilon)\big| \ge   (1-2\epsilon)\exp(H(\kappa)m - \epsilon m)\big|\FA(\sigma, \alpha:F,\epsilon)\big|$$
for all sufficiently large $m$. Thus,
\begin{eqnarray*}
H(\Sigma,\alpha \times \beta:F, 2\epsilon) \ge  H(\kappa) -\epsilon + H(\Sigma,\alpha:F,\epsilon).
\end{eqnarray*}
Let $\epsilon \to 0$ to obtain $H(\Sigma,\alpha \times \beta:F) \ge  H(\kappa)  + H(\Sigma,\alpha:F)$. This provides the lower bound and completes the proof.
\end{proof}

Theorem \ref{thm:mainapp} follows immediately from the above theorem, since it implies proposition \ref{prop:Bernoulli}, that $h(\Sigma, G,K^G,\kappa^G) = H(\kappa)$ whenever $H(\kappa) <\infty$. Of course, this uses theorem \ref{thm:K}, that $h(\Sigma,\cdot)$ is a measure-conjugacy invariant. To prove corollary \ref{cor:stepin}, we will need the following.

\begin{prop}
Let $G$ be a countable sofic Ornstein group. Let $(K,\kappa)$ be a standard Borel probability space with $H(\kappa) = +\infty$. Suppose there exists a generating partition $\alpha$ for $K^G$ such that $H(\alpha)<\infty$. Let $\Sigma$ be a sofic approximation to $G$. Then $h(\Sigma,G,K^G,\kappa^G)=-\infty$. 
\end{prop}

\begin{proof}
Let $(L,\lambda)$ be a probability space with $0<H(\lambda)<\infty$. By the previous proposition,
$$h\big(\Sigma, G, K^G \times L^G, \kappa^G \times \lambda^G\big) = h(\Sigma,G, K^G, \kappa^G) + H(\lambda).$$
There is a canonical measure-conjugacy between $(G, K^G \times L^G, \kappa^G \times \lambda^G)$ and $\big(G, (K \times L)^G, (\kappa \times \lambda)^G)$. Because $G$ is Ornstein and $H(\kappa \times \lambda)=+\infty$, this system is measurably conjugate to $(G,K^G,\kappa^G)$. Thus 
$$h\big(\Sigma, G, K^G \times L^G, \kappa^G \times \lambda^G\big) = h(\Sigma,G, K^G, \kappa^G).$$
From proposition \ref{prop:factor}, it follows that $h(\Sigma,G,K^G,\kappa^G) \le H(\alpha)<+\infty$. To see this, let $\beta$ be the trivial partition $\beta=\{X,\emptyset\}$. Since $H(\lambda)>0$ this implies $H(\Sigma,G,K^G,\kappa^G)=-\infty$ as claimed.
\end{proof}

Corollary \ref{cor:stepin} now follows from theorem \ref{thm:mainapp} and the proposition above. To prove theorem \ref{thm:generating} we will need:
\begin{cor}\label{cor:fbound}
Let $(G,X,\mu), (G,Y,\nu)$ be $G$-systems with finite-entropy generating partitions $\alpha, \beta$ respectively. If $(G,X,\mu)$ factors onto $(G,Y,\nu)$ and $\Sigma$ is a sofic approximation to $G$ then
$$h(\Sigma, G, Y,\nu) \ge h(\Sigma, G, X,\mu) - H(\alpha).$$
\end{cor}

\begin{proof}
Without loss of generality, we can assume that $X=Y$ and $\nu$ equals $\mu$ restricted to the $\sigma$-algebra generated by $\beta$ and the action of $G$. By proposition \ref{prop:factor}, for every finite $F \subset G$,
$$H(\Sigma,\beta:F) - H(\beta) \ge H(\Sigma, \alpha \vee \beta:F) - H(\alpha \vee \beta).$$
Since $H(\alpha \vee \beta) \le H(\alpha) + H(\beta)$ this implies 
$$H(\Sigma,\beta:F) \ge H(\Sigma,\alpha \vee\beta:F) - H(\alpha).$$
Take the infimum over all $F \subset G$ to obtain $h(\Sigma,\beta) \ge  h(\Sigma,\alpha \vee \beta)- H(\alpha)$. By theorem \ref{thm:K}, $h(\Sigma,\beta)=h(\Sigma,G,Y,\nu)$ and $h(\Sigma,\alpha \vee \beta)=h(\Sigma,G,X,\mu)$. So this proves it.
\end{proof}
We can now prove theorem \ref{thm:generating}.
\begin{proof}[Proof of theorem \ref{thm:generating}]
Let $G$ be a countable sofic group that contains a nonabelian free group. Let $(K,\kappa)$ be a standard Borel probability space with $H(\kappa)=+\infty$ and let $(L,\lambda)$ be a probability space with $0<H(\lambda)<\infty$. By theorem \ref{thm:weaklyisomorphic}, $(G,L^G,\lambda^G)$ factors onto $(G,K^G,\kappa^G)$. If the latter has a finite-entropy generating partition, the previous corollary implies that for any sofic approximation $\Sigma$ to $G$,
$$h(\Sigma, G,K^G,\kappa^G) \ge h(\Sigma, G, L^G, \lambda^G) - H(\lambda) = 0.$$
This contradicts the previous proposition. So $(G,K^G,\kappa^G)$ does not admit a finite-entropy generating partition.
\end{proof}

{\bf Acknowledgments}: I would like to thank Russ Lyons for suggesting that I think about the isomorphism problem for Bernoulli shifts over a nonabelian free group and for many useful conversations along the way. I'd also like to thank Benjy Weiss for asking whether the infinite entropy Bernoulli shift over a nonabelian free group could be finitely generated and for help with the history of entropy theory. I'd like to thank Charles Radin for asking whether every ensemble of circle packings of the hyperbolic plane can be approximated by finite packings. The answer lead to [Bo03] which eventually lead to the current work. I'm very grateful to Michael Hochman for an insightful conversation.

\end{document}